\documentclass[12pt,leqno]{amsart}

\usepackage{amsmath}
\usepackage{amssymb}
\usepackage{amsthm}
\usepackage{enumitem}
\usepackage{float}
\usepackage[T1]{fontenc}
\usepackage{geometry}
\usepackage{mathtools}
\usepackage{microtype}
\usepackage[table]{xcolor}
\usepackage[all,cmtip]{xy}

\def\acts{\ \rotatebox[origin=c]{-90}{$\circlearrowright$}\ }
\def\racts{\ \rotatebox[origin=c]{90}{$\circlearrowleft$}\ }

\usepackage[alphabetic]{amsrefs}
\usepackage[
    colorlinks,
    linkcolor=blue,
    anchorcolor=blue,
    citecolor=green
]{hyperref}
\usepackage[capitalise]{cleveref}

\theoremstyle{plain}
\newtheorem{theorem}{Theorem}[section]

\newtheorem{conjecture}[theorem]{Conjecture}
\newtheorem{corollary}[theorem]{Corollary}
\newtheorem{lemma}[theorem]{Lemma}
\newtheorem{proposition}[theorem]{Proposition}

\theoremstyle{definition}
\newtheorem{definition}[theorem]{Definition}

\newtheorem{remark}[theorem]{Remark}
\newtheorem{setup}[theorem]{}
\newtheorem*{nota-term}{Notation and Terminology}

\numberwithin{equation}{section}

\newcommand{\RomanNumeralCaps}[1]{\MakeUppercase{\romannumeral #1}}

\newcommand{\bk}{\mathbf{k}}

\DeclareMathOperator{\reg}{reg}
\DeclareMathOperator{\Exc}{Exc}
\DeclareMathOperator{\Gal}{Gal}
\DeclareMathOperator{\N}{N}
\DeclareMathOperator{\Nef}{Nef}
\DeclareMathOperator{\NS}{NS}
\DeclareMathOperator{\NE}{NE}
\DeclareMathOperator{\id}{id}
\DeclareMathOperator{\Frac}{Frac}
\DeclareMathOperator{\Pic}{Pic}

\DeclareMathOperator{\Spec}{Spec}
\DeclareMathOperator{\Supp}{Supp}

\begin{document}

\title[Surjective endomorphisms]
{Surjective endomorphisms of projective surfaces: the existence of infinitely many dense orbits}

\author{Jia Jia}
\address{National University of Singapore, Singapore 119076, Republic of Singapore}
\email{jia\_jia@u.nus.edu}

\author{Junyi Xie}
\address{BICMR, Peking University, Haidian District, Beijing 100871, China}
\email{xiejunyi@bicmr.pku.edu.cn}

\author{De-Qi Zhang}
\address{National University of Singapore, Singapore 119076, Republic of Singapore}
\email{matzdq@nus.edu.sg}


\begin{abstract}
	Let $f \colon X \to X$ be a surjective endomorphism of a normal projective surface.
	When $\deg f \geq 2$, applying an (iteration of) $f$-equivariant minimal model program (EMMP),
	we determine the geometric structure of $X$.
	Using this, we extend the second author's result to singular surfaces to the extent
	that either $X$ has an $f$-invariant non-constant rational function,
	or $f$ has infinitely many (disjoint) Zariski-dense forward orbits;
	this result is also extended to adelic topology (which is finer than Zariski topology).
\end{abstract}

\maketitle
\tableofcontents

\section{Introduction}

We work over an algebraically closed field $\bk$ of characteristic zero.
We first give a structure theorem for non-isomorphic surjective endomorphisms.

\begin{theorem}\label{thm:main}
	Let $f \colon X \to X$ be a non-isomorphic surjective endomorphism
	of a normal projective surface.
	Then $X$ has only log canonical (lc) singularities.
	If the canonical divisor $K_X$ is pseudo-effective, then Case~\ref{thm:main:qe} occurs.
	If $K_X$ is not pseudo-effective,
	replacing $f$ by an iteration,
	we may run an $f$-equivariant minimal model program (EMMP)
	\[
		X = X_1 \longrightarrow \cdots \longrightarrow X_j \longrightarrow \cdots \longrightarrow X_r \longrightarrow Y,
	\]
	contracting $K_{X_j}$-negative extremal rays,
	with $X_j \to X_{j+1}$ ($j < r$) being divisorial and $X_r \to Y$ being Fano contraction,
	such that one of the following cases occurs.
	\begin{enumerate}[label=(\arabic*),ref=1.1(\arabic*),wide=0pt,leftmargin=*]
		\item \label{thm:main:qe}
		      $f$ is quasi-{\'e}tale, i.e., {\'e}tale in codimension $1$;
		      there exists an $f$-equivariant quasi-{\'e}tale finite Galois cover
		      $\nu \colon V \to X$
		      as in \cref{thm:quasi-etale}.
		\item \label{thm:main:finite-ord}
		      $Y$ is a smooth projective curve of genus $g(Y) \geq 1$;
		      $f$ descends to an automorphism of finite order on the curve $Y$.
		\item \label{thm:main:smooth}
		      $Y$ is an elliptic curve;
		      $X \to Y$ is a $\mathbb{P}^1$-bundle.
		\item \label{thm:main:prdt}
		      $Y \cong \mathbb{P}^1$;
		      there is an $f|_{X_r}$-equivariant finite surjective morphism
		      $X_r \to Y \times \mathbb{P}^1$.
		\item \label{thm:main:cover}
		      $Y \cong \mathbb{P}^1$;
		      $f$ is polarized;
		      $K_X + S \sim_{\mathbb{Q}} 0$ for an $f^{-1}$-stable reduced divisor $S$;
		      there is an $f$-equivariant quasi-{\'e}tale finite Galois cover $\nu \colon V \to X$
		      as in \cref{thm:nak-7-1-1}\ref{thm:nak-7-1-1:4}.
		\item \label{thm:main:base-change}
		      $Y \cong \mathbb{P}^1$;
		      $f$ is polarized;
		      there exists an equivariant commutative diagram:
		      \[
			      \xymatrix{
			      \widetilde{f} \acts \widetilde{X} \ar@<2.3ex>[d]_{\widetilde{\pi}} \ar[r]^{\mu_{X}}	&	X \ar@<-2.3ex>[d]^{\pi} \racts f	\\
			      g_E \acts E \ar[r]												&	Y \racts g							\\
			      }
		      \]
		      here $E$ is an elliptic curve;
		      $\widetilde{X}$ is the normalisation of $X \times_Y E$;
		      $\widetilde{\pi}$ is a $\mathbb{P}^1$-bundle;
		      $\pi$ is a $\mathbb{P}^1$-fibration;
		      $\widetilde{f}$ and $g_E$ are finite surjective endomorphisms;
		      $\mu_X$ is quasi-{\'e}tale.
		\item \label{thm:main:fano-type}
		      $Y \cong \mathbb{P}^1$;
		      $X_r$ is of Fano type;
		      there is an $f|_{X_r}$-equivariant birational morphism $X_r \to \overline{X}$
		      to a klt Fano surface with Picard number $\rho(\overline{X}) = 1$;
		      every $X_j$ ($1 \leq j \leq r$) is a rational surface
		      whose singularities are klt (hence $\mathbb{Q}$-factorial).
		\item \label{thm:main:proj-cone}
		      $Y$ is a point;
		      $X_r$ is a projective cone over an elliptic curve $E$;
		      the normalisation $\Gamma$ of the graph of $X \dashrightarrow E$ is a $\mathbb{P}^1$-bundle,
		      and $f$ lifts to $f|_{\Gamma}$
		      such that $(f|_{\Gamma})^* |_{\N^1(\Gamma)} = \delta_f \id$.
		\item \label{thm:main:fano}
		      $Y$ is a point and hence $\rho(X_r)=1$;
		      $-K_{X_r}$ is ample;
		      every $X_j$ ($1 \leq j \leq r$) is a rational surface
		      whose singularities are lc and rational (hence $\mathbb{Q}$-factorial).
	\end{enumerate}
\end{theorem}

\begin{remark}
	By a result of Nakayama \cite{nakayama2020structure}*{Theorem~1.1},
	all the cases (1), (5) and (8) occur.
	Note that in case (8) above, the surface \(X_r\) is singular.
	Taking \(X\) to be the product of \(\mathbb{P}^1\) with a curve \(Y\) as in cases (2) -- (4) and (6),
	and a suitable non-isomorphic endomorphism \(f\) on \(X\),
	we see that cases (2) -- (4) and (6) are possible.
	A Hirzebruch surface
	\(\mathbb{F}_n \coloneqq \mathbb{P}(\mathcal{O}_{\mathbb{P}^1}\oplus\mathcal{O}_{\mathbb{P}^1}(n))\)
	with \(n\geq 1\) is an example of case (7) (cf.~\cite{nakayama2002ruled}*{Theorem~3}).

	Consider the weighted projective space \(X = \mathbb{P}(a_0,a_1,a_2)\)
	with coordinate functions \(x_0,x_1,x_2\)
	and the endomorphism \(f\colon X \to X\) defined by
	\[
		(x_0:x_1:x_2)\longmapsto (F_0(x_0,x_1,x_2):F_1(x_0,x_1,x_2):F_2(x_0,x_1,x_2))
	\]
	where \(F_j\)'s are homogeneous polynomials with respect to the weights
	such that \((\deg F_0: \deg F_1: \deg F_2) = (a_0:a_1:a_2)\)
	and \(F_j\)'s have no common zero.
	When \(\deg F_j \geq 2\) for some \(0\leq j \leq 2\),
	the pair \((X,f)\) gives an instance of case (9) above,
	for example, taking the singular surface \(X = \mathbb{P}(1,2,3)\)
	and \(f(x_0:x_1:x_2) = (x_0^2:x_1^2:x_2^2)\).
\end{remark}

The main ingredients of the proof for \cref{thm:main} are an equivariant minimal model program,
\cite{nakayama2020normal}*{Theorem A} and \cite{nakayama2020structure}*{Theorem 3.11}
(= \cref{thm:nak-7-1-1,thm:quasi-etale}),
and some analysis of Fano fibration in \cite{meng2019kawaguchi}*{Theorem 5.4}.
\cref{thm:main} implies:

\begin{corollary}\label{cor:main_compact}
	Let $f \colon X \to X$ be a non-isomorphic surjective endomorphism
	of a normal projective surface.
	Then, replacing $f$ by an iteration, one of the following holds.
	\begin{enumerate}[label=(\arabic*),ref=1.2(\arabic*),wide=0pt,leftmargin=*]
		\item \label{cor:main_compact:fib}
		      $X$ has an $f$-invariant non-constant rational function.
		\item \label{cor:main_compact:lifts}
		      $f$ lifts to an endomorphism on $X'$ via a generically finite surjective morphism,
		      where $X'$ is an abelian surface or a $\mathbb{P}^1$-bundle over an elliptic curve.
		\item \label{cor:main_compact:descends}
		      $f$ descends to an endomorphism on $\mathbb{P}^1 \times \mathbb{P}^1$
		      (commuting with both projections) via a generically finite surjective morphism.
		\item \label{cor:main_compact:pic-one}
		      $f$ descends to an endomorphism on a normal projective surface $\overline{X}$
		      via a birational morphism such that
		      \begin{enumerate}[label=(\alph*),ref=1.2(4\alph*)]
			      \item \label{cor:main_compact:a}
			            The Picard number $\rho(\overline{X}) = 1$,
			            so $f|_{\overline{X}}$ and hence $f$ are polarized;
			      \item \label{cor:main_compact:b}
			            The anti-canonical divisor $-K_{\overline{X}}$
			            is an ample $\mathbb{Q}$-Cartier divisor; and
			      \item \label{cor:main_compact:c}
			            $\overline{X}$ is a rational surface whose singularities
			            are log canonical (lc) and rational.
		      \end{enumerate}
	\end{enumerate}
\end{corollary}

Below is the motivation of our paper where \cref{conj:ZD}\ref{conj:ZD:extend}
is invariant under birational conjugation
which is stronger than the long-standing \cref{conj:ZD}\ref{conj:ZD:classical}.
In fact, \cref{conj:ZD}\ref{conj:ZD:classical} and \ref{conj:ZD:extend} are equivalent
modulo the Dynamical Mordell-Lang Conjecture as shown in \cite{xie2022existence}*{Proposition 2.6}
and also reminded by Professor Ghioca.

\begin{conjecture}\label{conj:ZD}
	Let $X$ be an (irreducible) projective variety over $\bk$
	and $f \colon X \dashrightarrow X$ a dominant rational self-map such that $\bk(X)^f = \bk$.
	Then:
	\begin{enumerate}[label=(\arabic*)]
		\item \label{conj:ZD:classical}
		      there is a point $x \in X(\bk)$ such that the (forward) orbit
		      $\mathcal{O}_f(x) \coloneqq \{f^s(x) \colon s \geq 0\}$ is well-defined,
		      i.e., $f$ is defined at $f^n(x)$ for any $n \geq 0$, and Zariski-dense in $X$.
		\item \label{conj:ZD:extend}
		      for every Zariski-dense open subset $U$ of $X$,
		      there exists a point $x\in X(\bk)$ whose orbit $O_f(x)$ under $f$ is well-defined,
		      contained in $U$ and Zariski-dense in $X$.
	\end{enumerate}
\end{conjecture}

To wit, a rational function $\psi \in \bk(X)$ is \emph{$f$-invariant}
if $f^*(\psi) \coloneqq \psi \circ f = \psi$.
Denote by $\bk(X)^f$ the field of $f$-invariant rational functions on $X$.
We have
\[
	\bk \subseteq \bk(X)^f \subseteq \bk(X).
\]

We will extend the classical \cref{conj:ZD} (Zariski-topology version)
to a stronger \cref{conj:AZD} (adelic-topology version).
We begin with:

\begin{setup}[\textbf{Adelic topology}]
	\label{ppty-of-adelic-topo}
	Assume that the transcendence degree of $\bk$ over $\mathbb{Q}$ is finite.
	In \cite{xie2022existence}*{\S~3},
	the second author has proposed the adelic topology on $X(\bk)$.
	The adelic topology has the following basic properties
	(cf.~\cite{xie2022existence}*{Proposition 3.16}).
	\begin{enumerate}[label=(\arabic*),wide=0pt,leftmargin=*]
		\item It is stronger than the Zariski topology.
		\item It is ${\mathsf{T}}_1$,
		      i.e., for any two distinct points $x, y \in X(\bk)$
		      there are adelic open subsets $U, V$ of $X(\bk)$
		      such that $x \in U, y \notin U$ and $y\in V, x \notin V$.
		\item Morphisms between algebraic varieties over $\bk$
		      are continuous for the adelic topology.
		\item Flat morphisms are open with respect to the adelic topology.
		\item \label{ppty-of-adelic-topo:irreducible}
		      The irreducible components of $X(\bk)$ in the Zariski topology
		      are the irreducible components of $X(\bk)$ in the adelic topology.
		\item Let $K$ be any subfield of $\bk$ which is finitely generated over $\mathbb{Q}$
		      and such that $X$ is defined over $K$ and $\overline{K} = \bk$.
		      Then the action
		      \[
			      \Gal(\bk/K)\times X(\bk)\to X(\bk)
		      \]
		      is continuous with respect to the adelic topology.
	\end{enumerate}
	When $X$ is irreducible, \ref{ppty-of-adelic-topo:irreducible} implies that
	the intersection of finitely many nonempty adelic open subsets of $X(\bk)$ is nonempty.
	So, if $\dim X\geq 1$, the adelic topology is not Hausdorff.
	In general, the adelic topology is strictly stronger than the Zariski topology.
\end{setup}

The adelic version of the Zariski-dense orbit conjecture was proposed in \cite{xie2022existence}.

\begin{conjecture}\label{conj:AZD}
	Assume that the transcendence degree of $\bk$ over $\mathbb{Q}$ is finite.
	Let $X$ be an irreducible variety over $\bk$ and
	$f \colon X \dashrightarrow X$ a dominant rational map.
	If $\bk(X)^f= \bk$, then there exists a nonempty adelic open subset $A \subseteq X(\bk)$
	such that the orbit of every point $x \in A$ is well-defined and Zariski-dense in $X$.
\end{conjecture}

\begin{definition}
	Let $X$ be an (irreducible) projective variety over $\bk$
	and $f \colon X \dashrightarrow X$ a dominant rational map.
	We say that a pair $(X,f)$ satisfies
	\begin{enumerate}
		\item \emph{ZD-property}, if \cref{conj:ZD}\ref{conj:ZD:classical} holds true;
		\item \emph{strong ZD-property}, if \cref{conj:ZD}\ref{conj:ZD:extend} holds true;
		\item \emph{AZD-property}, if \cref{conj:AZD} holds true; and
		\item \emph{SAZD-property}, if there is a nonempty adelic open subset $A$ of $X(\bk)$
		      such that for every point $x\in A$, its orbit $O_f(x)$ under $f$ is well-defined
		      and Zariski-dense in $X$.
	\end{enumerate}
\end{definition}

\begin{remark}\label{rem:AZD_to_ZD}
	\cref{conj:AZD} implies \cref{conj:ZD}.
	Precisely, we have:
	\begin{enumerate}
		\item SAZD-property implies AZD-property.
		\item \cref{conj:AZD} (adelic-topology version) is stronger
		      than the classical \cref{conj:ZD} (Zariski-topology version).
		      Indeed, even the hypothesis on $\bk$ in \cref{conj:AZD} does not cause any problem.
		      To be precise,
		      for every pair $(X,f)$ over $\bk$, there exists an algebraically closed subfield $K$
		      of $\bk$ whose transcendence degree over $\mathbb{Q}$ is finite
		      and such that $(X,f)$ is defined over $K$,
		      i.e., there exists a pair $(X_K,f_K)$ such that $(X,f)$ is its base change by $\bk$.
		      By \cite{xie2022existence}*{Corollary 3.29},
		      if $(X_K,f_K)$ satisfies AZD-property, then $(X,f)$ satisfies strong ZD-property.
	\end{enumerate}
\end{remark}

We will prove \cref{conj:AZD} (and hence \cref{conj:ZD}) for surjective endomorphisms of (possibly singular) projective surfaces, extending the smooth case in \cite{xie2022existence}.

\begin{theorem}\label{thm:AZD}
	Let $f \colon X \to X$ be a surjective endomorphism of a projective surface $X$ defined over the algebraically closed field $\bk$.
	Assume that $\bk$ has finite transcendence degree over $\mathbb{Q}$.
	Then \cref{conj:AZD} holds for $(X, f)$.
	Precisely, either $\bk(X)^f \neq \bk$;
	or there is a nonempty adelic open subset $A \subseteq X(\bk)$ such that the forward orbit $O_f(x)$ of every point $x \in A$ is Zariski-dense in $X$.
\end{theorem}

\begin{proposition}\label{pro:AZDinforbits}
	Assume that the transcendence degree of $\bk$ over $\mathbb{Q}$ is finite.
	Let $X$ be an irreducible variety over $\bk$ of positive dimension and $f \colon X \dashrightarrow X$ a dominant rational map.
	If $(X, f)$ has SAZD-property, then $f$ has infinitely many Zariski-dense orbits which are disjoint.
\end{proposition}

\cref{thm:AZD} and \cref{pro:AZDinforbits} imply the following (cf.~\cref{rem:AZD_to_ZD}).

\begin{theorem}\label{thm:ZD}
	Let $f \colon X \to X$ be a surjective endomorphism of a projective surface $X$ defined over the algebraically closed field $\bk$.
	Then \cref{conj:ZD} holds for $(X, f)$.
	Precisely, either $\bk(X)^f \neq \bk$;
	or for every Zariski-dense open subset $U$ of $X$, there is a point $x\in X(\bk)$ whose forward orbit $O_f(x)$ under $f$ is contained in $U$ and Zariski-dense in $X$.
	In the latter case, $f$ has infinitely many Zariski-dense orbits which are disjoint.
\end{theorem}

\par \noindent
\textbf{Ingredient of the proof.}
Due to the possible occurrence of Case (4) in \cref{cor:main_compact},
we cannot always reduce \cref{conj:AZD} to the smooth case.
Therefore, the following result for amplified endomorphisms is indispensable when proving \cref{thm:AZD}.
Since our $X$ might be singular, extra care is taken in the last section in proving it,
extending the smooth case in \cite{xie2022existence}.

\begin{theorem}\label{prop:int_amp}
	Let $X$ be a projective surface over $\bk$.
	Let $f \colon X \to X$ be an amplified endomorphism such that $\deg f > \delta_f$.
	Then the pair $(X,f)$ satisfies SAZD-property.

	In particular, if $f$ is int-amplified or polarized, the pair $(X,f)$ satisfies SAZD-property.
\end{theorem}

\begin{remark}
	Below are some histories of \cref{conj:ZD}\ref{conj:ZD:classical} (Zariski-topology version).
	It was proved by Amerik and Campana \cite{amerik2008fibrations}*{Theorem 4.1}
	when the field $\bk$ is uncountable.
	If $\bk$ is countable, it has been open even for the case of singular surfaces
	(with $f$ a well-defined morphism).
	Below are some confirmed cases.
	\begin{enumerate}[wide=0pt,leftmargin=*]
		\item In \cite{amerik2011existence}*{Corollary 9}, in arbitrary dimension,
		      Amerik proved the existence of non-preperiodic algebraic point
		      when $f$ is of infinite order.
		\item \Cref{conj:AZD} and hence \cref{conj:ZD} are proved
		      for $f = (f_1,\cdots,f_n) \colon (\mathbb{P}^1)^n \to (\mathbb{P}^1)^n$,
		      where the $f_i$'s are endomorphisms of $\mathbb{P}^1$,
		      in \cite{xie2022existence}*{Appendix B, joint work with Thomas Tucker};
		      see also \cite{bell2016dynamical}*{Theorem 14.3.4.2},
		      when $f_i$'s are not post-critically finite,
		      and Medvedev and Scanlon \cite{medvedev2014invariant}*{Theorem 7.16}
		      for endomorphism $f \colon \mathbb{A}^n \to \mathbb{A}^n$
		      with $f(x_1,\cdots,x_n) = (f_1(x_1),\cdots,f_n(x_n)),\, f_i(x_i) \in \bk[x_i]$.
		\item In \cite{xie2017existence}*{Theorem 1.1}, the second author has proved
		      \cref{conj:ZD}\ref{conj:ZD:classical} for dominant polynomial endomorphism
		      $f \colon \mathbb{A}^2 \to \mathbb{A}^2$.
		\item If $X$ is a (semi-) abelian variety and $f$ is a dominant self-map,
		      \cref{conj:ZD}\ref{conj:ZD:classical} has been proved in
		      \cite{ghioca2017density}*{Theorem 1.2} and \cite{ghioca2019density}*{Theorem 1.1};
		      in this abelian variety case, \cref{conj:AZD}
		      is proved in \cite{xie2022existence}*{Theorem 1.14}.
		\item When $X$ is an algebraic surface and $f$ is a birational self-map,
		      \cref{conj:AZD} (and hence \cref{conj:ZD}) has been proved in
		      \cite{xie2022existence}*{Corollary 3.31};
		      see also \cite{bell2014applications}*{Theorem 1.3} when $X$ is quasi-projective
		      over $\overline{\mathbb{Q}}$ and $f$ is an automorphism.
		\item \cref{conj:AZD} and hence \cref{conj:ZD} have been proved
		      when $X$ is a \emph{smooth} projective surface
		      and $f$ is a surjective endomorphism in \cite{xie2022existence}*{Theorem 1.15}.
	\end{enumerate}
\end{remark}

\par \noindent
{\bf Acknowledgements.}
The first and third authors are respectively supported by a President's Graduate Scholarship,
and the ARF A-8000020-00-00, from NUS.
The second author is partially supported by the project ``Fatou'' ANR-17-CE40-0002-01.
We would like to thank the referee for very careful reading, critical questions, and the valuable suggestions
to improve the paper.

\section{Preliminaries}

In this section, we collect together some definitions and preliminary results.

\begin{nota-term}
We use the following notation throughout the paper, with $X$ a projective variety.
\begin{table}[H]
	\begin{tabular}{p{5em}p{32em}}
		$\Pic X$      & the group of Cartier divisors of $X$ modulo linear equivalence $\sim$                                                                             \\
		$\Pic^0 X$    & the group of Cartier divisors of $X$ algebraically equivalent to 0                                                                                \\
		$\NS(X)$      & $\Pic X/\Pic^{0} X$, the N{\'e}ron-Severi group                                                                                                   \\
		$\N^1(X)$     & $\NS(X) \otimes_{\mathbb{Z}} \mathbb{R}$, the space of $\mathbb{R}$-Cartier divisors modulo numerical equivalence $\equiv$                        \\
		$\N_1(X)$     & the space of 1-cycles with coefficients in $\mathbb{R}$, modulo numerical equivalence $\equiv$                                                    \\
		$\Nef(X)$     & the cone of nef classes in $\N^1(X)$                                                                                                              \\
		$\NE(X)$      & the cone of pseudo-effective classes in $\N_1(X)$                                                                                                 \\
		$\kappa(X,D)$ & the Iitaka $D$-dimension of a $\mathbb{Q}$-Cartier divisor $D$	(cf.~\cite{iitaka1982algebraic}*{\S 10.1})                                         \\
		$\rho(X)$     & Picard number of $X$, which is $\dim_{\mathbb{R}} \N^1(X)$                                                                                        \\
		$q(X)$        & the irregularity of $X$, which is $h^1(X,\mathcal{O}_X) \coloneqq \dim H^1(X, \mathcal{O}_X)$                                                     \\
		$\Supp D$     & the support of an effective Weil divisor $D = \sum a_i D_i$ on $X$, which is $\cup_i D_i$, where $a_i >0$ and $D_i$'s are prime divisors          \\
		$f|_{Y}$      & the lifted (resp.~descended) endomorphism on \(Y\) of an endomorphism \(f\) on \(X\) via an equivariant morphism \(Y \to X\) (resp.~\(X \to Y\)).
	\end{tabular}
\end{table}
\end{nota-term}

An algebraic variety $X$ is called \emph{$\mathbb{Q}$-factorial} if every Weil $\mathbb{Q}$-divisor is $\mathbb{Q}$-Cartier.

An algebraic variety $X$ is said to have \emph{rational singularities} if $X$ is normal and if $R^i\pi_*\mathcal{O}_Y = 0$ ($i \geq 1$)
for one (and hence every) resolution
$\pi\colon Y \to X$ of singularities.

A \emph{pair} $(X,\Delta)$ consists of a normal variety $X$ and an effective Weil $\mathbb{R}$-divisor $\Delta=\sum b_i D_i$ such that $K_X + \Delta$ is $\mathbb{R}$-Cartier.
Let $\pi \colon Y\to X$ be a log resolution of the pair $(X,\Delta)$.
There are uniquely defined numbers $a(E_j,X,\Delta)$, called the \emph{discrepancy} of $E_j$ with respect to $(X,\Delta)$, such that
\[
	K_Y + \pi_*^{-1} \Delta = \pi^*(K_X + \Delta) + \sum_{E_j:\text{exceptional}} a(E_j,X,\Delta) E_j.
\]
We say that the pair $(X,\Delta)$ is \emph{Kawamata log terminal} (klt) (resp.~log canonical (lc))
if all $b_i < 1$ (resp.\ $\leq 1$), and, for one (and hence every) log resolution $\pi \colon Y \to X$ of $(X, \Delta)$,
we have $a(E_j, X,\Delta) > -1$ (resp.\ $\geq -1$) for every $\pi$-exceptional prime divisor $E_j$
(cf.~\cite{kollar1998birational}*{Definition 2.34, Corollary 2.32}).
We say that $X$ is klt or lc if so is $(X, 0)$.

Let $f \colon X \to X$ be a surjective endomorphism of a normal projective variety.
Then $f$ is finite and $f_* f^* = (\deg f) \id$ on $\Pic_{\mathbb{Q}} X$ (cf.~\cite{meng2019kawaguchi}*{Proposition 3.7}).

Let $\pi \colon X_1 \to X_2$ be a finite surjective morphism between normal varieties,
which automatically restricts to a finite morphism
$\pi' \coloneqq \pi|_{X_1'} \colon X_1' \to X_2'$ between smooth varieties
such that $X_i \setminus X_i'$ is a closed subvariety of $X_i$ of codimension $\ge 2$.
For a Weil divisor $D$ on $X_2$, we define $\pi^*D$ as the Zariski closure of $(\pi')^*(D|_{X_2'})$ in $X_1$.

Denote by $R_{\pi}$ the ramification divisor of $\pi$.
By definition, it is an effective Weil divisor and satisfies the \emph{ramification divisor formula:}
\begin{equation}\label{eq:ram-div-formula}
	K_{X_1} = \pi^{*}K_{X_2} + R_{\pi}.
\end{equation}
We say that $\pi$ is \emph{quasi-{\'e}tale} if it is {\'e}tale in codimension $1$, i.e., if $R_{\pi}=0$.

More generally, suppose $\pi^{-1}(D_2) = D_1$ for reduced effective Weil divisors $D_j \subset X_j$.
Then we have the \emph{logarithmic ramification divisor formula}
\begin{equation}\label{eq:log-ram-div-formula}
	K_{X_1} + D_1 = \pi^*(K_{X_2} + D_2) + R_\pi',
\end{equation}
where the \emph{log ramification divisor} $R_\pi'$ is an effective Weil divisor having no common irreducible component with $D_1$ (cf.~\cite{iitaka1982algebraic}*{\S 11.4}).
If $\pi^*D_2 = qD_1$, then $R_{\pi} = (q-1) D_1 + R_{\pi}'$.

Let $X$ be a normal projective surface which is $\mathbb{Q}$-factorial, and $C$ an irreducible curve on $X$.
Then $C$ is called a \emph{negative curve} if the self-intersection $C^2 < 0$.

A surjective morphism $f\colon X\to Y$ of algebraic varieties is a \emph{fibration} if $f$ has connected fibres.
Such an $f$ is a \emph{$\mathbb{P}^1$-fibration} if the general fibre of $f$ is isomorphic to $\mathbb{P}^1$.
A \emph{cross-section} of a fibration $f\colon X\to Y$ is an irreducible subvariety $C\subset X$ such that the restriction $f|_C\colon C\to Y$ is an isomorphism onto $Y$.

\begin{definition}\label{defn:fano-type}
	Let $X$ be a normal projective variety.
	Then $X$ is of \emph{Fano type} if there is an effective Weil $\mathbb{Q}$-divisor $\Delta$ such that $(X,\Delta)$ is klt and $-(K_X + \Delta)$ is ample.
\end{definition}

\begin{definition}\label{defn:dyn-pol}
	Let $f \colon X \to X$ be a surjective endomorphism of a projective variety.
	\begin{enumerate}[wide=0pt,leftmargin=*]
		\item The \emph{first dynamical degree} $\delta_f$ is the spectral radius of the endomorphism $f^*|_{\N^1(X)} \colon$ $\N^1(X) \to \N^1(X)$, i.e., the maximum of moduli of eigenvalues of $f^*|_{\N^1(X)}$.
		      It is known that $\deg f \leq \delta_f^{\dim X}$
		      (so $\deg f > 1$ implies $\delta_f > 1$).
		      For the case of the surface, see e.g.,~\cite{nakayama2020normal}*{Proposition 3.3}.
		      When $\dim X = 1$, we have $\delta_f = \deg f$.

		\item The map $f$ is \emph{$\delta$-polarized} if $f^*H \sim \delta H$ for some integer $\delta > 1$ and ample Cartier divisor $H$, or equivalently $f^*B \equiv \delta B$ for some rational number $\delta > 1$ and big $\mathbb{R}$-Cartier divisor $B$, or equivalently $f^*B \equiv \delta B$ for some $\delta > 1$ and big $\mathbb{Q}$-Cartier divisor $B$ (indeed such $\delta$ is in $\mathbb{Q}$) (cf.~\cite{meng2018building}*{Proposition 3.6}).
	\end{enumerate}
\end{definition}

For the following, see e.g.,~\cite{meng2019kawaguchi}*{Lemma~2.4},
\cite{meng2018building}*{Proposition2.9, Lemma~3.1, Corollary~3.12}
and \cite{zhang2010polarized}*{Lemma~2.2}.

\begin{lemma}\label{lem:dyn_pol}
	Let $f_i \colon X_i \to X_i$ ($i = 1, 2$) be surjective morphisms of projective varieties, $\pi \colon X_1 \dashrightarrow X_2$ a generically finite dominant rational map such that $\pi \circ f_1 = f_2 \circ \pi$.
	Then:
	\begin{enumerate}
		\item $\delta_{f_1} = \delta_{f_2}$.
		\item Suppose that $f_1^*B \equiv \delta B$ for some nef and big $\mathbb{R}$-Cartier divisor $B$ and $\delta > 0$.
		      Then $\deg f_1 = \delta^{\dim X_1}$.
		      In particular, if $f_1$ is $\delta$-polarized then $\deg f_1 = \delta^{\dim X_1}$.
		\item $f_1$ is $\delta$-polarized if and only if so is $f_2$.
		\item If $f_1$ is $\delta$-polarized,
		      then $f_1^*|_{\N^1(X)}$ is diagonalisable over \(\mathbb{C}\) and every eigenvalue of $f_1^*|_{\N^1(X)}$ has modulus $\delta$ ($> 1$).
		      In partially, $\delta_{f_1} = \delta$.
	\end{enumerate}
\end{lemma}

\begin{definition}
	Let $f \colon X \to X$ be a surjective endomorphism of a projective variety.
	\begin{enumerate}[wide=0pt,leftmargin=*]
		\item $f$ is \emph{amplified} if $f^*L - L = H$ for some Cartier divisor $L$ and ample divisor $H$.
		\item $f$ is \emph{int-amplified} if $f^*L - L = H$ for some ample Cartier divisors $L$ and $H$, or equivalently, if all the eigenvalues of $f^*|_{\N^1(X)}$ are of modulus greater than $1$ (cf.~\cite{meng2020building}*{Theorem 1.1}).
	\end{enumerate}
\end{definition}

\begin{remark}\label{rem:pol_to_amp}
	Let $f \colon X \to X$ be a surjective endomorphism of a projective variety.
	\begin{enumerate}[wide=0pt,leftmargin=*]
		\item If $f$ is polarized then $f$ is int-amplified.
		\item If $f$ is int-amplified and $\dim X \ge 2$, then $\deg f > \delta_f$ (cf.~\cite{meng2020building}*{Lemma 3.6}).
	\end{enumerate}
\end{remark}

The following results are frequently used.

\begin{lemma}(cf.~\cite{wahl1990characteristic}*{Theorems 2.8, 2.9}, \cite{nakayama2020normal}*{Theorem E}, \cite{boucksom2012volume}*{Theorem B})\label{lem:lc}
	Let $f \colon X \to X$ be a non-isomorphic surjective endomorphism of a normal variety of dimension two.
	Then we have:
	\begin{enumerate}
		\item $X$ has at worst lc singularities.
		      In particular, $K_X$ is $\mathbb{Q}$-Cartier.
		\item If $x_0 \in X$ is a closed point such that $f^{-1}(x_0) = x_0$ and $x_0 \in R_f$,
		      then $X$ has only klt singularity at $x_0$.
	\end{enumerate}
\end{lemma}

\begin{lemma}\label{lem:klt}
	Let $f$ be a surjective endomorphism of a normal projective surface $X$.
	Assume that there exists an $f$-equivariant fibration $\pi \colon X \to Y$ to a nonsingular projective curve
	such that $f$ descends to $g$ on $Y$ with $\deg g > 1$
	(this last part, $\deg g > 1$, automatically holds when $f$ is polarized by \cref{lem:dyn_pol}).
	Then $X$ has only klt singularities and hence is $\mathbb{Q}$-factorial.
\end{lemma}

\begin{proof}
	By \cref{lem:lc}, $X$ has at worst lc singularities.
	Assume $X$ is not klt (but is lc) at $x_0$.
	Then $f^{-1}(x_0) = x_0$ after iterating $f$ (cf.~\cite{broustet2014singularities}*{Lemma 2.10}).
	Also $x_0 \notin \Supp R_f$ by \cref{lem:lc}.
	Now $\pi \circ f = g \circ \pi$ implies $g^{-1}(\pi(x_0)) = \pi(x_0)$.
	Denote by $F_0 \coloneqq \pi^*(\pi(x_0))$ the fibre of $\pi$ passing through $x_0$.
	It follows that $f^{-1}(\Supp F_0) = \Supp F_0$ and thus $f^*F_0 = (\deg g) F_0$.
	But this yields $x_0 \in \Supp F_0 \subseteq \Supp R_f$, a contradiction.
\end{proof}

The general result below is a direct consequence of the cone theorem.

\begin{proposition}\label{prop:fibre-rational}
	Let $\pi \colon X \to Y$ be a morphism from a normal projective surface to a smooth projective curve with general fibres smooth rational curves.
	Then we have:
	\begin{enumerate}[label=(\arabic*)]
		\item \label{prop:fibre-rational:rat-sing}
		      $X$ has only rational singularities (and hence is $\mathbb{Q}$-factorial).
		      If $X$ has no negative curve, then $\Nef(X) = \NE(X)$.
		\item \label{prop:fibre-rational:fibre-irreducible}
		      Suppose that $\pi \colon X \to Y$ is a Fano contraction.
		      Then all the fibres are irreducible.
	\end{enumerate}
\end{proposition}

\begin{proof}
	The first part of \ref{prop:fibre-rational:rat-sing} is from \cite{nakayama2017variant}*{Proposition 2.33}.
	By ~\cite{zhang2016n}*{Lemma 2.3}, we have a natural embedding $\N^1(X) \subseteq \N_1(X)$.
	Since $X$ is $\mathbb{Q}$-factorial, we may identify $\N^1(X)$ and $\N_1(X)$.
	The second part follows from the same proof as \cite{kollar1996rational}*{\RomanNumeralCaps{2}, Lemma 4.12}.

	For \ref{prop:fibre-rational:fibre-irreducible}, suppose that $F_1, F_2$ are two distinct irreducible components of a fibre $F$.
	The $K_X$-negative extremal ray $R$ contracted by $\pi$ is $\mathbb{R}_{\geq 0} [F]$.
	Then $F_i \cdot F = 0$ implies $F_i \cdot R = 0$ and hence $F_i = \pi^*L_i$ for some $\mathbb{Q}$-Cartier divisors $L_i$ on $Y$ by the cone theorem \cite{fujino2011fundamental}*{Theorem 1.1(4) iii}.
	Since in our case $\Supp(L_i) = \pi(F_i) = \pi(F)$, it follows that $\Supp F_1 = \pi^{-1}(\pi(F)) = \Supp F_2$ as sets, a contradiction.
\end{proof}

Next, we deal with surface fibrations to curves of higher genus.

\begin{proposition}\label{prop:fib-to-elliptic}
	Let $f \colon X \to X$ be a surjective morphism of a normal projective surface.
	Suppose $\pi \colon X \to E$ is an $f$-equivariant fibration, with connected fibres, to an elliptic curve.
	We also assume that $X$ is $\mathbb{Q}$-factorial (this is the case when $\pi$ is a $\mathbb{P}^1$-fibration; cf.~\cref{prop:fibre-rational}).
	Let $g = f|_E$ and $\Sigma_1$ (resp.~$\Sigma_2$) be the set of points $e \in E$ such that the scheme-theoretic fibre $\pi^*(e)$ is reducible (resp.~nonreduced).
	Then $g^{-1}(\Sigma_i) = \Sigma_i$ for $i = 1, 2$.
\end{proposition}

\begin{proof}
	The result is clear when $f$ is an automorphism.
	Suppose $\deg f \geq 2$.
	Since general fibres of $\pi$ are smooth, the $\Sigma_i$'s are finite sets.
	It suffices to show $g^{-1}(\Sigma_i) \subseteq \Sigma_i$ for $i = 1, 2$.
	A fibre is reducible if and only if all its irreducible components are negative curves.
	After an iteration, $f^{-1}$ preserves negative curves (cf.~\cite{nakayama2002ruled}*{Lemma~9},
	or~\cite{meng2019kawaguchi}*{Lemma 4.3}),
	so $f^{-1}$ takes a negative curve to another one.
	Thus $f^{-1}$ takes reducible fibres to reducible fibres, and $g^{-1}(\Sigma_1) \subseteq \Sigma_1$.
	Since $E$ is an elliptic curve and hence $g \colon E \to E$ is {\'e}tale, the reducedness of $\pi^*(e)$ implies that of $\pi^*(g(e))$, so $g^{-1}(\Sigma_2) \subseteq \Sigma_2$ (cf.~\cite{cascini2020polarized}*{Lemma 7.3}).
\end{proof}

\begin{proposition}\label{prop:nonzero-irregularity}
	Suppose that $f \colon X \to X$ is a non-isomorphic surjective endomorphism of a normal projective surface.
	Let $\pi \colon X \to Y$ be a $\mathbb{P}^1$-fibration to a nonsingular projective curve with genus $g(Y) \geq 1$ such that $f$ descends to an endomorphism $h$ on $Y$.
	Then either $h$ is an automorphism of finite order; or $g(Y) = 1$ and $\pi$ is a $\mathbb{P}^1$-bundle.
\end{proposition}

\begin{proof}
	If the genus $g(Y) \geq 2$, then $Y$ is of general type, and the endomorphism $h$ is an automorphism of finite order.

	Suppose $g(Y) = 1$.
	Since $X$ has only lc singularities (cf.~\cref{lem:lc}) and $\deg f \geq 2$, we may run a relative EMMP over $Y$ (cf.~\cite{meng2019kawaguchi}*{Theorem 4.7})
	\[
		X = X_1 \to X_2 \to \cdots \to X_r \to Y',
	\]
	where $X_j \to X_{j+1}$ are divisorial contractions for $j < r$ and $X_r \to Y'$ is a Fano contraction.
	Since both $X \to Y$ and $X_r \to Y'$ have connected fibres,
	we have $Y' = Y$.
	Assume $r \geq 2$ and let $E$ be the exceptional divisor of $X_1 \to X_2$.
	Then $f^{-1}$ fixes $E$ as a set and thus $h^{-1}$ fixes $P$,
	where $P \coloneqq \pi(E)$ is a point in $Y$ (cf.~\cite{cascini2020polarized}*{Lemma 7.3}).
	Since $h$ is {\'e}tale, it must be an automorphism.
	Then $h$ is of finite order since it has at least one fixed point (cf.~\cite{hartshorne1977algebraic}*{\RomanNumeralCaps{4}, Corollary 4.7}).

	Assume that $h$ is not an automorphism of finite order.
	Then $r = 1$, i.e., $\pi$ is a Fano contraction.
	Keep the notation in \cref{prop:fib-to-elliptic}.
	\cref{prop:fibre-rational} implies $\Sigma_1 = \emptyset$.
	By \cref{prop:fib-to-elliptic}, one has $h^{-1}(\Sigma_2) = \Sigma_2$ and hence $\deg g \geq 2$ implies $\Sigma_2 = \emptyset$.
	So every fibre of $\pi$ is irreducible and reduced.
	Then $\pi$ is a smooth morphism (cf.~\cite{nakayama2017variant}*{Proposition 2.33}), and hence a $\mathbb{P}^1$-bundle (cf.~\cite{hartshorne1977algebraic}*{V.~Proposition~2.2}).
\end{proof}

For \cref{lem:proj-cone} below, $\mathfrak{Alb}(X)$ as defined in \cite{cascini2020polarized}*{Section 5}, is characterised by the rational map $\mathfrak{alb}_X \colon X \dashrightarrow \mathfrak{Alb}(X)$ as the one such that every rational map from $X$ to an abelian variety factors through it.

\begin{lemma}\label{lem:proj-cone}
	Let $X_r$ be a projective cone over an elliptic curve $E$.
	Let $X \to X_r$ be a birational morphism with $X$ a normal projective surface.
	Suppose $f \colon X \to X$ is a non-isomorphic surjective endomorphism and it descends to an endomorphism $f_r \colon X_r \to X_r$.
	Let $\Gamma$ be the normalisation of the graph of $\mathfrak{alb}_X \colon X \dashrightarrow E = \mathfrak{Alb}(X)$.
	Then $\Gamma \to E$ is a $\mathbb{P}^1$-bundle, $(f|_{\Gamma})^* = \delta_f \id$ and $\Gamma \to X_r$ is the contraction of a cross-section.
	Further, either $X = \Gamma$, or $X = X_r$.
\end{lemma}

\begin{proof}
	Note that $f$ descends to $E$ and lifts to $\Gamma$.
	Since $\rho(X_r) = 1$, $f_r$, and hence $f$, $f|_{\Gamma}$ and $f|_E$ are all polarized (cf.~\cref{lem:dyn_pol}).
	In particular, $f|_E$ is non-isomorphic.
	Applying \cref{prop:nonzero-irregularity}, $\Gamma \to E$ is a $\mathbb{P}^1$-bundle.
	Notice that $\rho(\Gamma) = 2$ and $\Gamma \to X_r$ is $f|_{\Gamma}$-equivariant.
	Let $C$ be the exceptional divisor for $\Gamma \to X_r$,
	which is a cross-section of $\Gamma \to E$.
	The exceptional divisor \(C\) is \((f|_{\Gamma})^{-1}\)-invariant,
	so $(f|_{\Gamma})^*C = \delta_f C$ (cf.~\cref{lem:dyn_pol}).
	Note that \(\rho(X_r) = 1\) and hence \((f|_r)^*|_{\N^1(X_r)} = \delta_f \id\) (cf.~\cref{lem:dyn_pol}).
	Since the pullback of $\N^1(X_r)$ is $(f|_{\Gamma})^*$-invariant
	and
	\(\N^1(\Gamma)\) is the direct sum of \(\mathbb{R}[C]\) and the pullback of \(\N^1(X_r)\),
	we conclude that $(f|_{\Gamma})^*$ is diagonalisable and thus $(f|_{\Gamma})^*|_{\N^1(\Gamma)} = \delta_f \id$.
	Since $2 = \rho(\Gamma) \ge \rho(X) \ge \rho(X_r) = 1$, Zariski main theorem and normality of $\Gamma$, $X$ and $X_r$ imply the last assertion.
\end{proof}

\cref{lem:K_X-pe-then-qe} below is known to Iitaka, Sommese, Fujimoto, Nakayama, $\cdots$.

\begin{lemma}\label{lem:K_X-pe-then-qe}
	Let $f$ be a non-isomorphic surjective endomorphism of a normal projective variety $X$ of dimension $n$.
	If $K_X$ is $\mathbb{Q}$-Cartier and pseudo-effective, then $R_f = 0$.
\end{lemma}

\begin{proof}
	Assume $R_f \neq 0$.
	The ramification divisor formula \eqref{eq:ram-div-formula} for $f^s$ is given by $K_X = (f^s)^*K_X + \sum_{i=0}^{s-1} (f^i)^*R_f$.
	Pick an ample Cartier divisor $H$ on $X$.
	Since $R_f$ is an integral Weil divisor and $K_X$ is pseudo-effective, we get a contradiction by letting $s \to +\infty$:
	\[
		K_X \cdot H^{n-1} = (f^s)^*K_X \cdot H^{n-1} + \sum_{i=0}^{s-1} (f^i)^*R_f \cdot H^{n-1} \geq s.\qedhere
	\]
\end{proof}

\Cref{thm:nak-7-1-1,thm:quasi-etale}
of Nakayama are crucial for the proof of \cref{thm:main}.

\begin{theorem}[cf.~\cite{nakayama2020normal}*{Theorem A}]\label{thm:nak-7-1-1}
	Let $X$ be a normal projective surface admitting a non-isomorphic surjective endomorphism $f$.
	Assume that $K_X + S \sim_{\mathbb{Q}} 0$ for an $f^{-1}$-stable reduced divisor $S$.
	Then $f$ is quasi-{\'e}tale outside $S$,
	and there exists a quasi-{\'e}tale finite Galois covering $\nu \colon V \to X$ such that
	$\nu \circ f_V = f^{\ell} \circ \nu$ for a non-isomorphic surjective endomorphism $f_V$ of $V$
	and a positive integer $\ell$,
	and that $V$ and $\nu$ satisfy exactly one of the following conditions.
	\begin{enumerate}[label=(\arabic*),wide=0pt,leftmargin=*]
		\item \label{thm:nak-7-1-1:3}
		      $V$ is an abelian surface and $S = 0$.
		\item \label{thm:nak-7-1-1:4}
		      $V$ is a $\mathbb{P}^1$-bundle over an elliptic curve
		      such that $\nu^* S$ is a disjoint union of two cross-sections.
		\item \label{thm:nak-7-1-1:5}
		      $V$ is a projective cone over an elliptic curve and $\nu^* S$ is a cross-section.
		\item \label{thm:nak-7-1-1:6}
		      $V$ is a toric surface and $\nu^* S$ is the boundary divisor.
	\end{enumerate}
\end{theorem}

\begin{theorem}[cf.~\cite{nakayama2020structure}*{Theorem 3.11}]\label{thm:quasi-etale}
	Let $X$ be a normal projective surface.
	Then $X$ admits a non-isomorphic quasi-{\'e}tale surjective endomorphism $f$
	if and only if there exists a quasi-{\'e}tale finite Galois covering $\nu\colon V\to X$
	satisfying one of the following conditions.
	\begin{enumerate}[label=(\arabic*),wide=0pt,leftmargin=*]
		\item \label{thm:quasi-etale:1}
		      $V$ is an abelian surface.
		\item \label{thm:quasi-etale:2}
		      $V \cong E \times T$ for an elliptic curve $E$ and a curve $T$ of genus at least two.
		\item \label{thm:quasi-etale:3}
		      $V \cong \mathbb{P}^1\times E$ for an elliptic curve $E$.
		\item \label{thm:quasi-etale:4}
		      $V$ is a $\mathbb{P}^1$-bundle over an elliptic curve associated
		      with an indecomposable locally free sheaf of rank two and degree zero.
	\end{enumerate}
	Moreover, $f^{\ell}$ lifts to $V$ for some positive integer $\ell$.
\end{theorem}

We need the following results of \cite{xie2022existence} in proving \cref{thm:AZD}.

\begin{proposition}[cf.~\cite{xie2022existence}*{Proposition 3.27}]\label{prop:AZD_bir}
	Let $X$ be an (irreducible) variety over $\bk$,
	and $f \colon X \dashrightarrow X$ a dominant rational map.
	Then the following statements are equivalent.
	\begin{enumerate}[wide=0pt,leftmargin=*]
		\item $(X, f)$ satisfies AZD-property (resp.~SAZD-property).
		\item $(X, f^m)$ satisfies AZD-property (resp.~SAZD-property) for some $m \geq 1$.
		\item There exists a pair $(Y, g)$ which is birational to the pair $(X, f)$,
		      and $(Y, g)$ satisfies AZD-property (resp.~SAZD-property).
	\end{enumerate}
\end{proposition}

\begin{lemma}[cf.~\cite{xie2022existence}*{Lemma 3.28}]\label{lem:AZD_gen}
	Let $X$ and $X'$ be (irreducible) varieties over $\bk$,
	$f \colon X \dashrightarrow X$ and $f' \colon X' \dashrightarrow X'$ dominant rational maps.
	Let $\pi \colon X' \dashrightarrow X$ be a generically finite dominant rational map
	such that $\pi \circ f' = f \circ \pi$.
	Then $(X', f')$ satisfies AZD-property (resp.~SAZD-property) if and only if
	$(X, f)$ satisfies AZD-property (resp.~SAZD-property).
\end{lemma}

\begin{proposition}[cf.~\cite{xie2022existence}*{Proposition 3.30}]\label{prop:SAZD_inf_curve}
	Let $X$ be an (irreducible) surface over $\bk$,
	and $f \colon X \dashrightarrow X$ a dominant rational map.
	Suppose the pair $(X, f)$ does not satisfy SAZD-property.
	Then there exists some $m \ge 1$
	such that there are infinitely many irreducible curves $C$ on $X$
	satisfying $f^m(C) \subseteq C$.
\end{proposition}

As a consequence of the above, we have the following.

\begin{corollary}\label{cor:AZD_aut} (cf.~\cite{xie2022existence}*{Corollary 3.31})
	Let $X$ be an (irreducible) projective surface over $\bk$, and $f \colon X \dashrightarrow X$ a birational map.
	Then the pair $(X, f)$ satisfies AZD-property.
\end{corollary}

\section{Proof of Theorem~\ref{thm:main} and Corollary~\ref{cor:main_compact}}

In this section, we prove \cref{thm:not-pseudo-eff} for the case of non-pseudo effective canonical divisor $K_X$.
\cref{thm:main} will follow from it and \cref{thm:nak-7-1-1} for the case with $K_X$ being pseudo-effective.

\begin{theorem}\label{thm:not-pseudo-eff}
	Let $f\colon X\to X$ be a non-isomorphic surjective endomorphism of a normal projective surface with $K_X$ not being pseudo-effective.
	Then $X$ has only lc singularities.
	Replacing $f$ by an iteration we may run an $f$-equivariant minimal model program
	\[
		\xymatrix{
		X \ar@{=}[r] \ar[d]_{f} & X_1 \ar[r]^{\pi_1} \ar[d]_{f_1} & \cdots \ar[r]^{\pi_{j-1}} & X_j \ar[r]^{\pi_j} \ar[d]_{f_j} & \cdots \ar[r]^{\pi_{r-1}} & X_r \ar[r]^{\pi} \ar[d]_{f_r} & Y \ar[d]^{g} \\
		X \ar@{=}[r] & X_1 \ar[r]_{\pi_1} & \cdots \ar[r]_{\pi_{j-1}} & X_j \ar[r]_{\pi_j} & \cdots \ar[r]_{\pi_{r-1}} & X_r \ar[r]_{\pi} & Y
		}
	\]
	contracting $K_{X_j}$-negative extremal rays, with $X_j \to X_{j+1}$ ($j < r$) being divisorial and $X_r \to Y$ being Fano contraction (hence every fibre of $\pi$ is irreducible).
	If $\dim Y = 1$, then $X_r$ has Picard number $\rho(X_r) = 2$ and all $X_j$ ($1 \leq j \leq r$) have only rational singularities (hence are $\mathbb{Q}$-factorial).
	Moreover, exactly one of the following cases occurs.
	\begin{enumerate}[label=(\arabic*),ref=3.1(\arabic*),wide=0pt,leftmargin=*]
		\item \label{thm:not-pseudo-eff:1} $Y$ is a smooth projective curve of genus $g(Y) \geq 2$,
		      and $g$ is an automorphism of finite order
		      (hence no $f_j$ is polarized by \cref{lem:dyn_pol},
		      and $r = 1$ by \cref{thm:mz-19}).
		\item \label{thm:not-pseudo-eff:2} $Y$ is an elliptic curve.
		      If $\pi$ has a nonreduced fibre, then $g$ is an automorphism of finite order;
		      otherwise $\pi \colon X_r \to Y$ is a $\mathbb{P}^1$-bundle.
		\item \label{thm:not-pseudo-eff:3}
		      $Y \cong \mathbb{P}^1$, and $f$ is not polarized (hence $r = 1$ by \cref{thm:mz-19}, so $X$ does not contain negative curves by \cref{lem:-K_X-0_2}, and $\Nef(X) = \NE(X)$ by \cref{prop:fibre-rational}).
		      \begin{enumerate}[label=(\alph*),ref=3.1(3\alph*)]
			      \item \label{thm:not-pseudo-eff:3a}
			            Either $\delta_f > \delta_g$;
			            or $\delta_f = \delta_g$, with $-K_X$ being ample or $R_f \neq 0$.
			            There exist a finite surjective morphism $\tau \colon X \to Y \times \mathbb{P}^1$
			            and a surjective endomorphism $h \colon \mathbb{P}^1 \to \mathbb{P}^1$ such that
			            $(g \times h) \circ \tau = \tau \circ f$.
			      \item \label{thm:not-pseudo-eff:3b}
			            $\delta_f = \delta_g$, $-K_X$ is nef but not ample, and $R_f = 0$.
			            There is an $f$-equivariant quasi-{\'e}tale finite Galois cover
			            $\nu \colon V \to X$ as in \cref{thm:quasi-etale}\ref{thm:quasi-etale:3}
			            or \ref{thm:quasi-etale:4}.
		      \end{enumerate}
		\item \label{thm:not-pseudo-eff:4}
		      $Y \cong \mathbb{P}^1$, and $f_r$ is polarized
		      (hence for all $1 \leq j \leq r$, $\delta_{f_j} = \delta_{f_r} = \delta_g$,
		      the $f_j$ is polarized by \cref{lem:dyn_pol},
		      and the $X_j$ has only klt singularities by \cref{lem:klt}).

		      If $X_r$ does not contain negative curves,
		      then $X_r$ satisfies \cref{prop:no-neg-curve}
		      and hence \ref{prop:no-neg-curve}\ref{prop:no-neg-curve:1}
		      or \ref{prop:no-neg-curve}\ref{prop:no-neg-curve:3} occurs.

		      Suppose that $X_r$ contains a negative curve $C$.
		      Then $\NE(X_r) = \langle [C], [F] \rangle$ with $F$ a general fibre of $\pi$;
		      $f_r^*C = \delta_f C$, so $R_{f_r} \geq (\delta_f - 1) C$;
		      also $\kappa(X_r,-K_{X_r}) = 2$;
		      further, one of the following is true.
		      \begin{enumerate}[start=3,label=(\alph*),ref=3.1{(4\alph*)}]
			      \item \label{thm:not-pseudo-eff:4c}
			            Either $-K_{X_r}$ is ample, or $C$ intersects $R_{f_r} - (\delta_f - 1)C$.
			            Then $X_r$ is of Fano type;
			            the contraction $\sigma \colon X_r \to \overline{X}$ of $C$ gives a klt Fano surface $\overline{X}$ with $\rho(\overline{X}) = 1$.
			            Moreover, $f_r$ descends to an endomorphism $\bar{f}$ on $\overline{X}$.
			      \item \label{thm:not-pseudo-eff:4d}
			            $-K_{X_r}$ is not ample, and $C$ does not meet $R_{f_r} - (\delta_f - 1)C$.
			            Then $C$ is a cross-section of $\pi$;
			            every irreducible component of $R_{f_r}$ dominates $Y$;
			            there exists an equivariant commutative diagram, where $E$ is an elliptic curve,
			            and $\widetilde{X}$ is the normalisation of (the main component of) $X \times_Y E$,
			            $\widetilde{\pi}$ is a $\mathbb{P}^1$-bundle;
			            $\widetilde{f}$ and $g_E$ are finite surjective endomorphisms;
			            $\mu_X$ is quasi-{\'e}tale, $\mu_Y$ is finite surjective;
			            $\overline{\pi}$ is the composition $X \to \cdots \to Y$.
			            \[
				            \xymatrix{
				            \widetilde{f} \acts \widetilde{X} \ar@<2.3ex>[d]_{\widetilde{\pi}} \ar[r]^{\mu_{X}}	&	X \ar@<-2.3ex>[d]^{\overline{\pi}} \racts f	\\
				            g_E \acts E \ar[r]_{\mu_Y}										&	Y \racts g									\\
				            }
			            \]
		      \end{enumerate}
		\item \label{thm:not-pseudo-eff:5}
		      $Y$ is a point, $-K_{X_r}$ is ample, $\rho(X_r) = 1$, and $f_r$ and hence all $f_j$ ($1 \leq j \leq r$) are polarized (cf.~\cref{lem:dyn_pol}).
		      Either $X_r$ is a projective cone over an elliptic curve, or a rational surface with only rational singularities.
	\end{enumerate}
\end{theorem}

\begin{remark}
	The $X_r$ has at most one negative curve when $\dim Y = 1$:
	if $C$ is a negative curve, then the class $[C]$ and fibre class are the only two extremal rays in $\NE(X_r)$.
\end{remark}

We need the following for the proof of \cref{thm:not-pseudo-eff}.

\begin{theorem}[cf.~\cite{meng2019kawaguchi}*{Theorem 5.4}]\label{thm:mz-19}
	Let $f$ be a non-isomorphic surjective endomorphism of a normal projective surface with $K_X$ not being pseudo-effective.
	Then $X$ has only lc singularities.
	Replacing $f$ by an iteration, we may run an $f$-equivariant MMP
	\[
		X = X_1 \to \cdots \to X_j \to \cdots \to X_r \to Y,
	\]
	contracting $K_{X_j}$-negative extremal rays, with $X_j \to X_{j+1}$ ($j < r$) being divisorial
	and $X_r \to Y$ being Fano contraction, such that one of the following cases occurs (with $f_j = f|_{X_j}$).
	\begin{enumerate}[label=(\arabic*),wide=0pt,leftmargin=*]
		\item \label{thm:mz-19:1}
		      $\dim Y = 0$, so $\rho(X_r) = 1$, also $f_r$ and hence all $f_j$ ($1 \leq j \leq r$) are polarized.
		\item \label{thm:mz-19:2}
		      $\dim Y = 1$, $f_r$ and hence all $f_j$ ($1 \leq j \leq r$) are polarized;
		      $\rho(X_r) = 2$ and $\delta_f = \delta_{f|Y}$.
		\item \label{thm:mz-19:3}
		      $\dim Y = 1$ and $f_r$ is not polarized.
		      Further $r = 1$, $\rho(X) = 2$ and one of the following cases occurs.
		      \begin{enumerate}[label={(\alph*)},ref={(3\alph*)}]
			      \item \label{thm:mz-19:3a}
			            $\delta_f = \delta_{f|Y}$.
			      \item \label{thm:mz-19:3b}
			            $\delta_f > \delta_{f|Y}$;
			            there exist a finite surjective morphism $\tau \colon X \to Y\times\mathbb{P}^1$
			            and a surjective endomorphism $h \colon \mathbb{P}^1 \to \mathbb{P}^1$
			            such that $(f|_Y \times h) \circ \tau = \tau \circ f$.
		      \end{enumerate}
	\end{enumerate}
\end{theorem}

\begin{remark}
	Let $f \colon X \to X$ be a surjective endomorphism of a normal projective variety, and $\pi \colon X \to Y$ an $f$-equivariant fibration, with connected fibres, to a smooth projective curve.
	Let $F$ be a general fibre of $\pi$.
	Then $f^* F \equiv \deg (f|_Y) F = (\delta_{f|Y}) F$.
	In particular, one of the eigenvalues of $f^*|_{\N^1(X)}$ is a positive integer.
\end{remark}

We will apply the \cref{lem:-K_X-pe} -- \cref{prop:no-neg-curve} below to $X = X_r$ in \cref{thm:not-pseudo-eff}.

\begin{lemma}\label{lem:-K_X-pe}
	Let $f \colon X \to X$ be a non-isomorphic surjective endomorphism of a normal projective surface.
	Let $\pi \colon X \to Y$ be a Fano contraction such that $f$ descends to an endomorphism $g$ on $Y$.
	Assume either $\dim Y = 0$;
	or $\dim Y = 1$ and $\delta_g=\deg g > 1$.
	Then $-K_X$ is pseudo-effective but not numerically trivial.
\end{lemma}

\begin{proof}
	First, by \cref{lem:lc}, $X$ has only lc singularities and $K_X$ is $\mathbb{Q}$-Cartier.
	If $\dim Y = 0$, then the Picard number $\rho(X) = 1$;
	hence $-K_X$ is ample, whence $K_X$ is not pseudo-effective.

	Suppose $\dim Y = 1$ and $\deg g > 1$.
	Then $\rho(X) = \rho(Y) + 1 = 2$.
	Write $\NE(X) = \langle [C],[F] \rangle$, where $F$ is a general fibre of $\pi$,
	and $[C]$ is another extremal divisor class.
	Since $K_X \cdot F = \deg (K_X|_F) = \deg K_F = -2$, we may assume $-K_X \equiv C + bF$ with $C \cdot F = 2$.
	Suppose to the contrary that $b = -b_1 < 0$.
	Write $f^*C \equiv \delta C$ in $\NE(X)$.
	We calculate the ramification divisor:
	\begin{align}
		R_f = K_X - f^*K_X & \equiv (\delta - 1) C - b_1(\delta_g - 1) F,\nonumber \\
		(\delta - 1) C     & \equiv R_f + b_1 (\delta_g - 1) F.	\label{eq:-K_X-pe}
	\end{align}
	Now we have reached a contradiction from \eqref{eq:-K_X-pe}, because $R_f$ is effective, $b_1(\delta_g - 1) > 0$ (by the assumption) and $[F]$ and $[C]$ are two distinct extremal divisor classes.
\end{proof}

\begin{lemma}\label{lem:-K_X-0_2}
	Let $f \colon X \to X$ be a non-isomorphic surjective endomorphism of a normal projective surface, and $\pi \colon X \to Y$ an $f$-equivariant Fano contraction with $\dim Y = 1$.
	Let $C$ be a negative curve on $X$.
	Then we have:
	\begin{enumerate}[label={(\alph*)}]
		\item \label{lem:-K_X-0_2:pol}
		      $f$ is polarized.

		\item \label{lem:-K_X-0_2:kappa}
		      If further, $Y \cong \mathbb{P}^1$, then $\kappa(X,-K_X) = 0$ or $2$.
	\end{enumerate}
\end{lemma}

\begin{proof}
	Note that $X$ is $\mathbb{Q}$-factorial by \cref{prop:fibre-rational}.
	By assumption, $\NE(X) = \langle [C], [F] \rangle$ where $F$ is a general fibre of $\pi$.
	Note that $f^*C \equiv \delta C$ for some $\delta > 0$ and hence $f^*C = \delta C$, since $C^2 < 0$.
	Then
	\begin{align*}
		(\deg f) C \cdot F & = f^*C \cdot f^*F = (\delta \,\delta_{f|Y}) C \cdot F, \\
		(\deg f) C^2       & = (f^*C)^2 = \delta^2 C^2
	\end{align*}
	imply that $\delta = \delta_{f|Y}$ ($= \deg f|_Y$) and $f$ is $\delta$-polarized.

	Now suppose $Y \cong \mathbb{P}^1$.
	Then $q(X) = 0$.
	Hence, by \cref{lem:-K_X-pe} and its proof, $-K_X$ is numerically and hence $\mathbb{Q}$-linearly equivalent to some effective divisor.
	Thus $\kappa(X,-K_X) \geq 0$.
	If $-K_X$ is not big, then its class lies in the boundary of $\NE(X)$.
	Since $-K_X \cdot F = 2$, $-K_X \sim_{\mathbb{Q}} aC$ for some $a>0$.
	Accordingly, $\kappa(X,-K_X) = \kappa(X,C) = 0$.
\end{proof}

Next, we consider the case when the $\mathbb{Q}$-factorial $X$ has no negative curve.

\begin{lemma}\label{lem:neg-curve}
	Let $f \colon X\to X$ be a non-isomorphic surjective endomorphism of a normal projective surface.
	Suppose that $\pi \colon X \to \mathbb{P}^1$ is an $f$-equivariant Fano contraction such that $\delta_f = \delta_{f|{\mathbb{P}^1}}$
	(hence $X$ is a rational surface, and has only klt singularities;~cf.~\cref{lem:klt}).
	Write $\Nef(X) = \langle [D],[F] \rangle$ with $F$ a general fibre of $\pi$, and let $f^*D \equiv \lambda D$.
	Then:
	\begin{enumerate}[label={(\alph*)}]
		\item \label{lem:neg-curve:ev-int}
		      $\lambda$ is a positive integer.
		\item \label{lem:neg-curve:D-big}
		      If $D^2 > 0$ then $f$ is polarized and $X$ has at least one negative curve.
		\item \label{lem:neg-curve:D-pencil}
		      Assume that $D^2 = 0$ (hence $D$ is a $\mathbb{Q}$-divisor after replacing it by a positive multiple) and $\kappa(X,D) > 0$.
		      Then $D$ is a semi-ample $\mathbb{Q}$-Cartier divisor and $\kappa(X,D) = 1$.
	\end{enumerate}
\end{lemma}

\begin{proof}
	By \cref{lem:klt}, $X$ is $\mathbb{Q}$-factorial.
	Since $f^*F \sim_{\mathbb{Q}} \delta_{f|{\mathbb{P}^1}} F = \delta_f F$, $0 < (\deg f) D \cdot F = f^*D \cdot f^*F = (\lambda \, \delta_f) D \cdot F$, so $\lambda \in \mathbb{Q}_{> 0}$ is an algebraic integer and hence an integer.


	For \ref{lem:neg-curve:D-big}, if $D^2>0$, then $(\lambda \, \delta_f) D^2 = (\deg f) D^2 = (f^*D)^2 = \lambda^2 D^2$ implies $\lambda = \delta_f$ and $f$ is polarized.
	Assume $X$ has no negative curve.
	Then $\Nef(X) = \NE(X)$ by \cref{prop:fibre-rational}.
	Thus $\mathbb{R}_{\geq 0}[D]$ is extremal in $\Nef(X) = \NE(X)$ and $D^2 = 0$, a contradiction.

	For \ref{lem:neg-curve:D-pencil},
	we have $\kappa(X,D) = 1$ since $D^2 = 0$.
	We may assume $h^0(X,D) \geq 2$, and write $|D| = |M| + F_D$ as the moving part and fixed part.
	Pick two general $D_1, D_2 \in |M|$, so $D_1$ and $D_2$ have no common components.
	Since $D$ and $M$ are nef,
	$0 = D^2 \ge D \cdot M \ge M^2 = D_1 \cdot D_2 \ge 0$.
	Thus $M \cdot F_D = 0$ and $D_1 \cdot D_2 = 0$, hence the $D_i$ are semi-ample.
	Now $0 = D^2 = (M + F_D)^2$ implies $F_D^2 = 0$ ($= M^2 = M \cdot F_D$).
	Hodge index theorem (cf.~\cite{nakayama2017variant}*{Lemma~pp.~302}) implies that $F_D$ is numerically
	and hence $\mathbb{Q}$-linearly equivalent to a multiple of the movable $M$.
	Replacing $D$ by a multiple, we may assume $F_D = 0$ and hence $D = M$ is semi-ample.
\end{proof}

\begin{proposition}\label{prop:no-neg-curve}
	With the same assumptions as in \cref{lem:neg-curve}, assume $X$ has no negative curve.
	Then $\Nef(X) = \NE(X) = \langle [D], [F] \rangle$.
	Also, one of the following cases occurs.
	\begin{enumerate}[label=(\arabic*),wide=0pt,leftmargin=*]
		\item \label{prop:no-neg-curve:1}
		      There exist a finite surjective morphism $\tau \colon X \to \mathbb{P}^1 \times \mathbb{P}^1$ and a surjective endomorphism $h \colon \mathbb{P}^1 \to \mathbb{P}^1$ such that $(f|_{\mathbb{P}^1} \times h)\circ \tau = \tau \circ f$.
		\item \label{prop:no-neg-curve:2}
		      $f$ is quasi-{\'e}tale, but non-polarized.
		      There is an $f^{\ell}$-equivariant quasi-{\'e}tale finite Galois cover $\nu \colon V \to X$ for an integer $\ell > 0$ as in \cref{thm:quasi-etale}\ref{thm:quasi-etale:3} or \ref{thm:quasi-etale:4}.
		\item \label{prop:no-neg-curve:3}
		      $f$ is polarized.
		      The (nef non-big) divisor $D$ can be chosen to be irreducible, reduced and $f^{-1}$-stable such that $K_X + D \sim_{\mathbb{Q}} 0$.
		      There is an $f^{\ell}$-equivariant quasi-{\'e}tale finite Galois cover $\nu \colon V \to X$ for an integer $\ell > 0$ as in \cref{thm:nak-7-1-1}\ref{thm:nak-7-1-1:4}.
	\end{enumerate}
\end{proposition}

\begin{proof}
	By \cref{lem:neg-curve}, $D^2 = 0$,
	so $D$ is nef and non-big.
	By \cref{prop:fibre-rational}, $\NE(X) = \Nef(X) = \langle [D],[F] \rangle$.
	In particular, $[D]$ generates an extremal ray in $\NE(X)$.
	Then, by \cref{lem:-K_X-pe}, $-K_X \in \NE(X) = \Nef(X)$ is a nef divisor.
	Note that $-K_X \cdot F = 2$.

	\par \vskip 1pc \noindent
	{\bf Case A.} $-K_X \cdot D > 0$, i.e., $-K_X$ is ample.
	By the cone theorem (cf.~\cite{kollar1998birational}*{Theorem 3.7}), $D$ may be chosen to be a rational curve.
	Since $(aD - K_X) \cdot D > 0$ and $(aD - K_X) \cdot F > 0$ for $a > 0$, Kleiman's ampleness criterion (cf.~\cite{kollar1998birational}*{Theorem 1.18}) implies that $aD - K_X$ is ample.
	Hence $D$ is semi-ample by the basepoint-free theorem (cf.~\cite{kollar1998birational}*{Theorem 3.3}),
	so it gives a fibration $\phi \colon X \to Z$ with connected fibres and normal $Z$.
	We have $\dim Z = 1$, since $D$ is not big.
	Notice that $F \cdot D > 0$, $F$ dominates $Z$ and hence $Z \cong \mathbb{P}^1$.

	Let $C$ be an irreducible curve on $X$.
	Then $C$ is in a fibre of $\phi$ if and only if $C \cdot D = 0$.
	Using the projection formula we obtain $f_*C \cdot D = C \cdot f^*D = \lambda (C \cdot D)$.
	Consequently, $C$ is in a fibre if and only if so is $f(C)$.
	Since the fibration $\phi$ has connected fibres, there exists a surjective endomorphism $h \colon Z \to Z$ such that $h \circ \phi = \phi \circ f$ by the rigidity lemma (cf.~\cite{debarre2001higher}*{Lemma 1.15}).
	The two distinct fibrations $\pi$ and $\phi$ induce a surjective morphism $\tau \colon X \to \mathbb{P}^1\times Z = \mathbb{P}^1\times \mathbb{P}^1$ such that $(f|_{\mathbb{P}^1} \times h) \circ \tau = \tau \circ f$.
	It is finite because $\rho(X) = 2 = \rho (\mathbb{P}^1 \times \mathbb{P}^1)$.
	Then Case~\ref{prop:no-neg-curve:1} occurs.

	\par \vskip 1pc \noindent
	\textbf{Case B.}
	$-K_X \cdot D = 0$, i.e., $-K_X$ is nef but not ample.
	Since $D$ is also nef, Hodge index theorem (cf.~\cite{nakayama2017variant}*{Lemma~pp.~302}) implies (the class of) $-K_X \in \mathbb{R}_{\geq 0}[D]$.
	Hence $f^*K_X \sim_{\mathbb{Q}} \lambda K_X$, because $q(X) = 0$.

	If $R_f = 0$ (which is impossible when $f$ is polarized, for otherwise $K_X = f^*K_X \equiv 0$ but $-K_X \cdot F = 2$), then Case~\ref{prop:no-neg-curve:2} occurs by \cref{thm:quasi-etale}.
	Indeed, only \cref{thm:quasi-etale}\ref{thm:quasi-etale:3} - \ref{thm:quasi-etale:4} occur since $-K_X$ is nef and not numerically trivial.

	Assume $0 \neq R_f$ ($= K_X - f^*K_X \in \mathbb{R}_{\geq 0} [D]$)
	and write $R_f = \sum a_i D_i$ where $a_i \in \mathbb{Z}_{> 0}$ and $D_i$ are irreducible components, automatically with $[D_i] \in \mathbb{R}_{\geq 0} [D]$.
	If $R_f$ is reducible, then $q(X) = 0$ implies $D_1 \sim_{\mathbb{Q}} tD_2$ for some rational number $t > 0$;
	or if $D_1$ is not $f^{-1}$-stable, then $f^*D_1 \sim_{\mathbb{Q}} \lambda D_1$.
	In either case, $D_1$ is $\mathbb{Q}$-Cartier, $D_1^2 = 0$ and $\kappa(X,D_1) > 0$.
	Hence by \cref{lem:neg-curve} we may retake $D \coloneqq D_1$ and assume it is a semi-ample $\mathbb{Q}$-Cartier divisor and $\kappa(X,D) = 1$.
	Using the semi-ample divisor $D$, Case~\ref{prop:no-neg-curve:1} occurs as argued in Case A.

	Now we may retake $D \coloneqq D_1$ (hence $D$ is a $\mathbb{Q}$-divisor) and assume $\Supp R_f = D$ is irreducible (and reduced) and $f^{-1}$-stable.
	\par \vskip 1pc \noindent
	{\bf Case B-1.} $f$ is not polarized.
	If $\lambda = 1$ we reach a contradiction: $0 \leq R_f = (K_X - f^*K_X) \sim_{\mathbb{Q}} 0$ and hence $R_f = 0$.
	Thus $1 < \lambda < \delta_f$.
	By \cite{meng2019kawaguchi}*{Proof of Step 5 in Theorem 5.2} we obtain $\kappa(X,D) > 0$.
	Arguing as above then Case~\ref{prop:no-neg-curve:1} occurs.

	\par \vskip 1pc \noindent
	{\bf Case B-2.} $f$ is polarized.
	We shall show that Case~\ref{prop:no-neg-curve:3} occurs.
	Now the log ramification divisor formula \eqref{eq:log-ram-div-formula} becomes $K_X + D = f^*(K_X + D)$.
	So the eigenvector (of $f^*$) $K_X + D \sim_{\mathbb{Q}} 0$ since $f$ is polarized and $q(X) = 0$.
	Applying \cref{thm:nak-7-1-1}, only Cases~\ref{thm:nak-7-1-1}\ref{thm:nak-7-1-1:4} - \ref{thm:nak-7-1-1:6} may occur in our situation.
	Note that $X$ has klt singularities by \cref{lem:klt}, so does $V$ and hence Case~\ref{thm:nak-7-1-1}\ref{thm:nak-7-1-1:5} cannot occur.
	For Case~\ref{thm:nak-7-1-1}\ref{thm:nak-7-1-1:6}, $\nu^* D$ would be the boundary divisor of a toric surface, so a big divisor.
	But this violates the fact that $\kappa(V,\nu^*D) = \kappa(X,D) < 2$.
	So only Case~\ref{thm:nak-7-1-1}\ref{thm:nak-7-1-1:4} is possible,
	i.e., Case~\ref{prop:no-neg-curve:3} occurs.
\end{proof}

Now we can prove the main result of this section.

\begin{proof}[Proof of \cref{thm:not-pseudo-eff}]
	The assertions in the first paragraph follow from \cref{lem:lc}, \cref{thm:mz-19} and \cref{prop:fibre-rational}.

	We now apply \cref{thm:mz-19}.
	If Case~\ref{thm:mz-19}\ref{thm:mz-19:1} occurs, i.e., if $\dim Y = 0$ then \ref{thm:not-pseudo-eff:5} occurs.
	Indeed, the last assertion there has been proved in \cite{broustet2017remarks}*{p.~578}.
	If $\dim Y = 1$ with $g(Y) \geq 1$,
	then $q(X_r) = g(Y) \geq 1$ and \ref{thm:not-pseudo-eff:1} -- \ref{thm:not-pseudo-eff:2} occur by \cref{prop:nonzero-irregularity}.

	Suppose that \ref{thm:mz-19}\ref{thm:mz-19:3} occurs with $Y \cong \mathbb{P}^1$.
	Then $f$ is not polarized and $r = 1$.
	\cref{lem:-K_X-0_2} implies $X$ does not contain negative curves.
	By \cref{prop:no-neg-curve} and its proof, Case~\ref{thm:mz-19}\ref{thm:mz-19:3a} leads to \ref{thm:not-pseudo-eff:3}:
	if $-K_X$ is ample or $R_f \neq 0$ then \ref{thm:not-pseudo-eff:3a} occurs;
	if $-K_X$ is not ample and $R_f = 0$ then \ref{thm:not-pseudo-eff:3b} occurs.
	Clearly, Case~\ref{thm:mz-19}\ref{thm:mz-19:3b} implies \ref{thm:not-pseudo-eff:3a}.

	Now suppose that \ref{thm:mz-19}\ref{thm:mz-19:2} occurs with $Y \cong \mathbb{P}^1$.
	Thus $f_r$ is polarized, and $X_r$ has only klt singularities by \cref{lem:klt}.
	The second paragraph of \ref{thm:not-pseudo-eff:4} follows from Case~\ref{thm:mz-19}\ref{thm:mz-19:2} and \cref{prop:no-neg-curve} when $X_r$ has no negative curve.

	We still have to consider Case~\ref{thm:mz-19}\ref{thm:mz-19:2} with the extra conditions that $Y \cong \mathbb{P}^1$, $X_r$ has only klt singularities and contains a negative curve $C$.
	By \cref{lem:-K_X-0_2}, $f_r$ is $\delta_f$-polarized (cf.~\cref{lem:dyn_pol}) and $\kappa(X_r,-K_{X_r}) = 0, 2$.
	Let $F$ be a general fibre of $\pi \colon X_r \to Y$.
	Since both $[C]$ and $[F]$ are extremal classes of $\NE(X_r)$ and $\rho(X_r) = 2$, $\NE(X_r) = \langle [C], [F] \rangle$.
	Note that $-K_{X_r} \cdot F = 2$, $f_r^* C = \delta_f C$.
	The log ramification divisor formula \eqref{eq:log-ram-div-formula} for $f_r$ is
	\begin{equation}\label{eq:log-ram-f_r}
		K_{X_r} + C = f_r^*(K_{X_r} + C) + R_{f_r}'
	\end{equation}
	with $R_{f_r} = (\delta_f - 1)C + R_{f_r}'$.

	Consider the case $\kappa(X_r,-K_{X_r}) = 0$.
	We will reach a contradiction.
	By the proof of \cref{lem:-K_X-0_2}, this happens if and only if $-K_{X_r} \sim_{\mathbb{Q}} aC$ for some $a > 0$.
	Then \eqref{eq:log-ram-f_r} gives
	\begin{equation}\label{eq:R_f'-C}
		R_{f_r}' \sim_{\mathbb{Q}} (a - 1) (\delta_f - 1) C.
	\end{equation}
	Since $R_{f_r}'$ is effective and has no common component with $C$, \eqref{eq:R_f'-C} gives $a = 1$ and $R_{f_r}' = 0$.
	The eigenvector (of $f_r^*$) $K_{X_r} + C \sim_{\mathbb{Q}} 0$ since $f_r$ is polarized and $q(X_r) = 0$.
	Now we apply \cref{thm:nak-7-1-1}.
	Only Case~\ref{thm:nak-7-1-1}\ref{thm:nak-7-1-1:4} may occur in our situation (cf.~Proof of Case B-2 in \cref{prop:no-neg-curve}).
	But for Case~\ref{thm:nak-7-1-1}\ref{thm:nak-7-1-1:4} we have $0 = K_V^2 = (\deg \nu) K_{X_r}^2 = (\deg \nu) (-C)^2 < 0$, a contradiction.
	So the case $\kappa(X_r,-K_{X_r}) = 0$ will not occur.

	Thus we may assume $\kappa(X_r,-K_{X_r}) = 2$.
	If $-K_{X_r}$ is ample, i.e., $-K_{X_r} \cdot C > 0$, then $X_r$ is a klt Fano surface.
	The $K_{X_r}$-negative extremal contraction $\sigma \colon X_r \to \overline{X}$ of $C$ gives a klt Fano surface $\overline{X}$ with $\rho(\overline{X}) = 1$
	(cf.~\cite{kollar1998birational}*{Corollary 3.43(1)}).
	Note that $f_r$ descends to an endomorphism $\bar{f} \colon \overline{X} \to \overline{X}$ since $f^{-1}$ fixes $C$ as a set.
	Thus Case~\ref{thm:not-pseudo-eff:4c} occurs.

	Next, we may also assume that $-K_{X_r}$ is not ample, i.e., $-K_{X_r} \cdot C \leq 0$.
	Let $-K_{X_r} = P + aC$ be the Zariski decomposition (cf.~\cite{sakai1984weil}*{Corollary~(7.5)}) with $P$ nef and big.
	We claim that $P \cdot C = 0$.
	Indeed, if $a > 0$, then the `negative part' $C$ has $P \cdot C = 0$;
	if $a = 0$, then $0 \leq P \cdot C = -K_{X_r} \cdot C \leq 0$, and hence $P \cdot C = 0$.

	Now we may assume that things like $X$, $f$, $C$ etc. are defined over a field \(K\)
	which is finitely generated over \(\mathbb{Q}\).
	Embedding \(K\) into \(\mathbb{C}\),
	we may assume that the base field is \(\mathbb{C}\) and we may use some techniques from complex-analytic geometry in the following.

	Let $\sigma \colon X_r \to \overline{X}$ be the contraction of the negative curve $C$
	to a point $\overline{x}$ in the normal Moishezon surface $\overline{X}$
	(cf.~\cite{sakai1984weil}*{Theorem~(1.2)}).
	Then $f_r$ descends to $\bar{f} \colon \overline{X} \to \overline{X}$ with $\bar{f}^{-1}(\bar x) = \bar{x}$.
	By \cref{lem:lc}, $\overline{X}$ has only lc singularities and $K_{\overline X}$ is $\mathbb{Q}$-Cartier.
	We have $P = \sigma^*(-K_{\overline X})$ since both sides are perpendicular to (and hence the same modulo) the negative curve $C$.
	Hence $\sigma^*K_{\overline X} = -P = K_{X_r} + aC$.

	Suppose that $C$ intersects $R_{f_r}'$.
	We shall show that Case~\ref{thm:not-pseudo-eff:4c} occurs.
	Indeed, then $\overline{x} \in \Supp R_{\bar f}$ and hence $(\overline{X},\overline{x})$ and also $\overline{X}$ have only klt singularities by \cref{lem:lc}, so $\overline{X}$ is $\mathbb{Q}$-factorial.
	For a Moishezon surface, being $\mathbb{Q}$-factorial implies the projectivity
	(cf.~\cite{fujino2021minimalb}*{Lemma~4.1}).
	Then $-K_{\overline X}$ is ample, and $\overline{X}$ is a klt Fano surface with $\rho(\overline{X}) = 1$.
	Since $\overline{X}$ has only klt singularities, $\sigma^*K_{\overline X} = K_{X_r} + aC$ implies that the pair $(X_r,aC)$ is klt.
	For $b \in \mathbb{Q}_{> 0}$ with $0 < b - a \ll 1$, the pair $(X_r,bC)$ is still klt (cf.~\cite{kollar1998birational}*{Corollary 2.35(2)});
	the divisor
	\[
		-(K_{X_r} + bC) = -(K_{X_r} + aC) - (b - a)C = P - (b - a)C
	\]
	has positive intersections with the two generators $[C]$ and $[F]$ of $\NE(X_r)$ and hence is ample.
	Thus $X_r$ is of Fano type.
	So Case~\ref{thm:not-pseudo-eff:4c} occurs.

	Suppose now that $C$ does not intersect $R_{f_r}'$ even if raise $f_r$ to some power.
	We shall show that Case~\ref{thm:not-pseudo-eff:4d} occurs.
	Note that $R_{f_r}' \neq 0$, otherwise by \eqref{eq:log-ram-f_r} the eigenvector (of the polarized $f_r^*$) $K_{X_r} + C \sim_{\mathbb{Q}} 0$, reaching a contradiction: $2 = \kappa(X_r,-K_{X_r}) = \kappa(X_r,C) = 0$.

	Note that $F$ is nef but not numerically trivial, so $C \cdot F > 0$.
	It follows that $C$ dominates $Y$, i.e., $\pi(C) = Y$.
	Since $R_{f_r}'$ is disjoint with $C$ by assumption and every fibre of $X_r \to Y$ is irreducible by \cref{prop:fibre-rational}, every irreducible component of $R_{f_r}'$ (and hence of $R_{f_r}$) is horizontal and hence dominates $Y$.
	Then Case~\ref{thm:not-pseudo-eff:4d} is true for $(X_r,f_r)$ by \cite{matsuzawa2019kawaguchi}*{Theorem 4.4} and also for $(X,f)$ by \cite{matsuzawa2019kawaguchi}*{Lemma 4.12}.
	In fact, it follows from \cref{prop:nonzero-irregularity} that $\widetilde{\pi}$ is a $\mathbb{P}^1$-bundle, since $f$ is polarized and hence $g_E$ cannot be an automorphism.
	The equality
	\begin{align*}
		-2 + F \cdot C & = F \cdot(K_{X_r} + C) = F \cdot (f_r^*(K_{X_r} + C) + R_{f_r}') \\
		               & = \delta_{f} (-2 + F \cdot C) + F \cdot R_{f_r}',
	\end{align*}
	and $F \cdot R_{f_r}' > 0$ imply that $C$ is a cross-section of $\pi$.
\end{proof}

\begin{lemma}\label{lem:rat-sing}
	Let $X$ be a normal projective surface with lc singularities.
	Let $\varphi \colon X \to X'$ be a composition of extremal contractions.
	Suppose that either $X'$ is a rational surface and $H^i(X',\mathcal{O}_{X'}) = 0$ for $i = 1, 2$;
	or $X' \cong \mathbb{P}^1$.
	Then $X$ is a rational surface whose singularities are rational.
\end{lemma}

\begin{proof}
	The first assertion is clear.
	For the second, we prove by induction on the number of contractions in $\varphi$.
	We may assume $\varphi$ itself is the contraction of a $K_X$-negative extremal ray.
	Then $-K_X$ is $\varphi$-ample.
	By the relative Kodaira vanishing theorem (cf.~\cite{fujino2011fundamental}*{Theorem 8.1}), $R^i\varphi_*\mathcal{O}_{X} = 0$ for every $i > 0$.
	Then the Leray spectral sequence
	\[
		E_2^{p,q} = H^p(X,R^q\varphi_*\mathcal{O}_{X}) \Rightarrow H^{p+q}(X',\mathcal{O}_{X'})
	\]
	degenerates and $H^i(X,\mathcal{O}_{X}) = H^i(X',\mathcal{O}_{X'}) = 0$ for $i = 1, 2$ in both cases.
	Since $X$ is a rational surface, it has rational singularities by \cite{broustet2017remarks}*{Lemma 2.7}.
\end{proof}

\begin{proof}[Proof of \cref{thm:main}]
	By \cref{lem:lc}, $X$ has lc singularities and $K_X$ is $\mathbb{Q}$-Cartier.
	If $K_X$ is pseudo-effective, then \cref{lem:K_X-pe-then-qe} implies that $f$ is quasi-{\'e}tale.
	Applying \cref{thm:quasi-etale},
	only Cases~\ref{thm:quasi-etale}\ref{thm:quasi-etale:1} - \ref{thm:quasi-etale:2}
	are possible since $K_X$ is pseudo-effective and $\nu$ is quasi-{\'e}tale.
	Therefore Case~\ref{thm:main:qe} holds.
	If $K_X$ is not pseudo-effective, then we can apply \cref{thm:not-pseudo-eff}.
	We have the correspondence below.
	\begin{table}[H]
		\begin{tabular}{>{\color{blue}}c|>{\color{blue}}c}
			\hline
			\ref{thm:main:qe}          & \ref{thm:quasi-etale}, \ref{thm:not-pseudo-eff:3b}                                                              \\ \hline
			\ref{thm:main:finite-ord}  & \ref{thm:not-pseudo-eff:1} -- \ref{thm:not-pseudo-eff:2}, \ref{prop:nonzero-irregularity}                       \\ \hline
			\ref{thm:main:smooth}      & \ref{thm:not-pseudo-eff:2}, \ref{prop:nonzero-irregularity}                                                     \\ \hline
			\ref{thm:main:prdt}        & \ref{thm:not-pseudo-eff:3a}, \ref{thm:not-pseudo-eff:4} (with \ref{prop:no-neg-curve}\ref{prop:no-neg-curve:1}) \\ \hline
			\ref{thm:main:cover}       & \ref{thm:not-pseudo-eff:4} (with \ref{prop:no-neg-curve}\ref{prop:no-neg-curve:3})                              \\ \hline
			\ref{thm:main:base-change} & \ref{thm:not-pseudo-eff:4d}                                                                                     \\ \hline
			\ref{thm:main:fano-type}   & \ref{thm:not-pseudo-eff:4c}                                                                                     \\ \hline
			\ref{thm:main:proj-cone}   & \ref{thm:not-pseudo-eff:5}, \ref{lem:proj-cone}                                                                 \\ \hline
			\ref{thm:main:fano}        & \ref{thm:not-pseudo-eff:5}                                                                                      \\ \hline
		\end{tabular}
	\end{table}

	Set $f_r \coloneqq f|_{X_r}$.
	We first show Case~\ref{thm:not-pseudo-eff:4} (with \ref{prop:no-neg-curve}\ref{prop:no-neg-curve:3}) implies Case~\ref{thm:main:cover}.
	By \cref{prop:no-neg-curve}\ref{prop:no-neg-curve:3}, $X_r$ is a rational surface; $f_r$ is polarized;
	there is an $f_r^{-1}$-stable reduced divisor $D \ne 0$ on $X_r$ such that $f_r$ is quasi-{\'e}tale outside $D$.
	Let $\pi \colon X \to \cdots \to X_r$ be the composition of the divisorial contractions, $D' \coloneqq \pi^{-1}_*(D)$ the proper transform of $D$ and $\Exc(\pi)$ the exceptional divisor of $\pi$.
	Set $S \coloneqq D' + \Exc(\pi)$.
	Then $S$ is reduced and $f^{-1}$-stable since $\pi$ is $f$-equivariant.
	Since the divisor $D$ on $X_r$ is not big and $\pi_*(S) = D$, the divisor $S$ is also non-big on $X$.
	By the definition of $S$, $f$ is quasi-{\'e}tale outside $S$.
	Thus $K_X + S = f^*(K_X + S)$,
	by the log ramification divisor formula \eqref{eq:log-ram-div-formula}.
	Since $f$ is polarized (cf.~\cref{lem:dyn_pol}) and $q(X) = 0$, the eigenvector (of $f^*$) $K_X + S \sim_{\mathbb{Q}} 0$.
	Note that $X$ has only klt singularities by \cref{lem:klt}.
	Now we can apply \cref{thm:nak-7-1-1} to say Case~\ref{thm:main:cover} occurs (cf.~Proof of Case B-2 in \cref{prop:no-neg-curve}).

	It remains to prove the last assertions of Cases~\ref{thm:main:fano-type} and \ref{thm:main:fano}, respectively.
	\cref{lem:klt} implies that $X_j$ has only klt singularities for $1 \leq j \leq r$ in Case~\ref{thm:main:fano-type}.
	In Case~\ref{thm:main:fano}, we have $H^i(X_r,\mathcal{O}_{X_r}) = 0$ for $i = 1, 2$ by \cite{fujino2011fundamental}*{Theorem 8.1}.
	These $X_j$'s in Case~\ref{thm:main:fano} have only rational singularities by \cref{lem:rat-sing}.
\end{proof}

\begin{proof}[Proof of \cref{cor:main_compact}]
	We apply \cref{thm:main}.
	Case~\ref{thm:main:qe} implies either Case~\ref{cor:main_compact:lifts} occurs;
	or $V = E \times T$ for an elliptic curve $E$ and a curve $T$ of genus $g(T) \ge 2$ and then $V$ has an $(f|_V)^m$-invariant non-constant rational function for some $m \ge 1$, hence $X$ has an $f$-invariant non-constant rational function by \cite{xie2022existence}*{Lemma 2.1}, and thus Case~\ref{cor:main_compact:fib} occurs.
	Case~\ref{thm:main:finite-ord} leads to Case~\ref{cor:main_compact:fib}.
	Cases~\ref{thm:main:smooth}, \ref{thm:main:cover} - \ref{thm:main:base-change} and \ref{thm:main:proj-cone} satisfy Case~\ref{cor:main_compact:lifts}.
	Case~\ref{thm:main:prdt} implies Case~\ref{cor:main_compact:descends} (replacing $f$ by $f^2$).
	Cases~\ref{thm:main:fano-type} and \ref{thm:main:fano} satisfy the conditions \ref{cor:main_compact:a} -- \ref{cor:main_compact:c} (cf.~\cref{lem:dyn_pol}) and thus Case~\ref{cor:main_compact:pic-one}.
\end{proof}

\section{Amplified endomorphisms; Proofs of Theorems~\ref{thm:AZD} and \ref{prop:int_amp} and Proposition~\ref{pro:AZDinforbits}}

In this section we assume that the transcendence degree of $\bk$ over $\mathbb{Q}$ is finite.
Till \cref{prop_gd_pt}, we fix an (irreducible) projective surface $X$ over $\bk$, and a surjective endomorphism $f \colon X \to X$.

The key is to prove \cref{prop:int_amp} which is essential in proving \cref{thm:AZD}.

\begin{definition}
	Let $o \in X(\bk)$ be a smooth fixed point of $f$, and $\lambda_1,\lambda_2$ the eigenvalues of the tangent map $df_o \coloneqq df|_o \colon T_{X,o} \to T_{X,o}$.
	The smooth point $o \in X(\bk)$ is said to be a \emph{repelling} fixed point of $f$ with respect to a norm $|\cdot|$ of $\bk$, if $|\lambda_i|>1$ for $i=1,2$.
	The smooth point $o \in X(\bk)$ is said to be a \emph{good} fixed point of $f$, if $df_o$ is invertible and one of the following conditions holds:
	\begin{enumerate}[label=(\arabic*)]
		\item
		      $\lambda_1$ and $\lambda_2$ are multiplicatively independent;
		\item \label{defn:good-pt:embd}
		      There exist a prime $p$ and an embedding $\tau \colon \bk \hookrightarrow \mathbb{C}_p$ such that
		      \[
			      |\tau(\lambda_1) + \tau(\lambda_2)|\leq1 \text{ and } |\tau(\lambda_1)||\tau(\lambda_2)|<1
		      \]
		      where $|\cdot|$ is the $p$-adic norm on $\mathbb{C}_p$.
	\end{enumerate}
\end{definition}

\begin{remark}
	Note that Condition~\ref{defn:good-pt:embd} just means that both $|\tau(\lambda_1)|$ and $|\tau(\lambda_2)|$ are at most one and $|\tau(\lambda_i)| < 1$ for $i = 1$ or $2$.
\end{remark}

\begin{definition}
	We say that $f$ has \emph{R-property} if there exist a smooth fixed point $o$ of $f$ and an embedding $\sigma \colon \bk \hookrightarrow \mathbb{C}$ such that both $|\sigma(\lambda_1)|$ and $|\sigma(\lambda_2)|$ are strictly greater than $1$, where $\lambda_1,\lambda_2$ are the eigenvalues of the tangent map $df_o \colon T_{X,o} \to T_{X,o}$.
\end{definition}

\textit{From now on till \cref{prop_gd_pt}, we assume that $f \colon X \to X$ is an amplified endomorphism, i.e., $f^*L\otimes L^{-1}$ is an ample line bundle for some line bundle $L$ on $X$.}

We will show the existence of a good fixed point of $f$.
Let $R$ be a finitely generated $\overline{\mathbb{Q}}$-sub-algebra of $\bk$, such that $\bk$ is the algebraic closure of $\Frac R$, and $X, f, L$ are defined over $\Frac R$.
There is a variety $X_{\Frac R}$ over $\Frac R$ and an endomorphism $f_{\Frac R} \colon X_{\Frac R}\to X_{\Frac R}$, such that $X = X_{\Frac R}\times_{\Spec {\Frac R} }\Spec \bk$ and $f = f_{\Frac R} \times_{\Spec {\Frac R}} \id$.

After shrinking $W \coloneqq \Spec R$, we may assume that $W$ is smooth, there exists a projective $R$-scheme $\pi \colon X_R\to W$ whose generic fibre is $X_{\Frac R}$, $f_{\Frac R}$ extends to a finite endomorphism $f_R$ on $X_R$ and there exists a line bundle $L_R$ on $X_R$ such that $f_R^*L_R\otimes L_R^{-1}$ is $\pi$-ample.
For every point $t\in W(\overline{\mathbb{Q}})$, denote by $X_t$ the special fibre $X_R\times_{W}\Spec \overline{\mathbb{Q}}$ of $X_R$ over $t$.
Let $L_t, f_t$ be the restrictions of $L_R, f_R$ on $X_t$.
Let $X_R^{\reg}$ be the smooth locus of $X_R$.
After shrinking $W$, we assume that $X_t$ is irreducible for every $t \in W(\overline{\mathbb{Q}})$ and $X_t \cap X_R^{\reg} \neq \emptyset$.

\begin{lemma}\label{lemliftgoodf}
	Assume that there exists some $t\in W(\overline{\mathbb{Q}})$ such that $f_t$ has a good fixed point in $X_t\cap X_R^{\reg}$.
	Then $f$ has a good fixed point in $X$.
\end{lemma}

\begin{proof}
	This lemma is shown by the proof of \cite{xie2022existence}*{Lemma 6.6}.
\end{proof}

\begin{proposition}\label{prop_gd_pt}
	Let $f\colon X\to X$ be an amplified endomorphism of a projective surface $X$ over $\bk$.
	Let $o$ be a smooth fixed point of $f$ such that $df_o$ is invertible.
	Let $C$ be an irreducible curve in $X$ passing through $o$ such that $f(C)=C$, and every branch of $C$ at $o$ is invariant under $f$.
	Denote by $\pi_C \colon \overline{C} \to C$ the normalisation of $C$ and $f|_{\overline C} \colon \overline{C} \to \overline{C}$ the endomorphism induced by $f|_C$.
	Let $q \in \pi_C^{-1}(o)$ and set $\mu \coloneqq d(f|_{\overline{C}})|_q$.
	Assume that there exists an embedding $\alpha \colon \bk \hookrightarrow \mathbb{C}$ such that $0<|\alpha(\mu)|<1$.
	Then there exists some $n \geq 1$ such that $f^n$ has a good fixed point in $X$.
\end{proposition}

\begin{proof}
	This proof is a small modification of the proof of \cite{xie2022existence}*{Lemma 6.7}.

	After enlarging $R$, we may assume $o,C,q$ are defined over $\Frac R$ and $\mu \in R$.
	After shrinking $W$, we may assume there is an irreducible subscheme $C_R$ of $X_R$ whose generic fibre is $C$ and a section $o_R \in X_R(R)$ of $\pi \colon X_R \to W$ whose
	$\bk$-extension is $o$.
	For every point $t \in W(\overline{\mathbb{Q}})$, denote by $C_t$ and $o_t$ the specializations of $C_R$ and $o_R$.
	Because $o \in X^{\reg}_R$, after shrinking $W$, we may assume that $o_t\in X^{\reg}_R$ and $C_t$ is irreducible for every $t\in W(\overline{\mathbb{Q}})$.
	There is a projective morphism $\pi_{C_R} \colon \overline{C}_R\to C_R$ over $R$ whose generic fibre is $\pi_{C}$ and an $R$-point $q_R\in \overline{C}_R(R)$, whose generic fibre is $q$.
	After shrinking $W$, we may assume that for all $t \in W(\overline{\mathbb{Q}})$, the specialization $\pi_{C_t}\colon \overline{C}_t\to C_t$ of $\pi_{C_R}$ is the normalisation of $C_t$.

	The embedding $\alpha \colon R\subseteq \bk\hookrightarrow \mathbb{C}$ defines a point $\eta \in W(\mathbb{C})$.
	We view $\mu$ as a function on $W(\mathbb{C})$.
	We have $|\mu(\eta)| = |\alpha(\mu)|\in (0,1)$.
	There exists a Euclidean open neighborhood $U$ of $\eta$, such that $|\mu(\cdot)|\in (0,1)$ on $U$.
	Picking $t \in U \cap W(\overline{\mathbb{Q}})$, we have $0 < |\mu(t)| < 1$.
	By \cref{lemliftgoodf}, we only need to prove that there exists some $n \geq 1$ such that $f_t^n$ has a good fixed point in $X_t$.
	Thus we have reduced to the case $\bk = \overline{\mathbb{Q}}$.

	\textit{Now we may assume that $\bk=\overline{\mathbb{Q}}$, the surface $X$ and the map $f$ are defined over a number field $K$, and there exist a variety $X_K$ over $K$ and an endomorphism $f_K \colon X_K \to X_K$, such that $X = X_K \times_{\Spec K} \Spec \bk$ and $f = f_K\times_{\Spec K} \id$.}

	Let $O_K$ be the ring of integers of $K$.
	There exists a projective $O_K$-scheme $X_{O_K}$ which is flat over $\Spec O_K$ whose generic fibre is $X_K$.
	Denote by $\pi_{O_K} \colon X_{O_K}\to \Spec O_K$ the structure morphism.
	The endomorphism $f_{K}$ on the generic fibre extends to a rational self-map $f_{O_K}$ on $X_{O_K}$.
	Denote by $o_{O_K}$ the Zariski closure of $o$ in $X_{O_K}$.
	Since $o$ is defined over $K$, $o_{O_K}$ is a section of $\pi_{O_K}$.

	Denote by $\pi^{O_K}_{\mathbb{Z}} \colon \Spec {O_K} \to \Spec \mathbb{Z}$ the morphism induced by the inclusion $\mathbb{Z} \hookrightarrow {O_K}$.
	Let $X_{\mathbb{Z}}$ be the $\mathbb{Z}$-scheme which is the same as $X_{O_K}$ as an absolute scheme with the structure morphism $\pi_{\mathbb{Z}} \coloneqq \pi^{O_K}_{\mathbb{Z}}\circ\pi_{O_K} \colon X_{O_K}\to \Spec \mathbb{Z}$.
	Then $X_{\mathbb{Z}}$ is a projective $\mathbb{Z}$-scheme.
	Denote by $f_{\mathbb{Z}} \colon X_{\mathbb{Z}}\dashrightarrow X_{\mathbb{Z}}$ the rational self-map induced by $f_{O_K}$.

	Since $f_{\mathbb{Z}}$ is regular on the generic fibre, there exists a finite set $B(f,\mathbb{Z})$ of primes such that $f_{\mathbb{Z}}$ is regular on $\pi_{\mathbb{Z}}^{-1}(\Spec \mathbb{Z} \setminus B)$.
	Set $B(f, {O_K}) \coloneqq (\pi^{O_K}_{\mathbb{Z}})^{-1}(B(f,\mathbb{Z}))$.
	Set $W \coloneqq \Spec O_K \setminus B(f,{O_K})$, $X_W \coloneqq \pi_{O_K}^{-1}(\Spec O_K \setminus B(f,{O_K}))$.
	Then $f_{O_K}$ is regular on $X_W$.
	Set $o_W \coloneqq o_{O_K}\cap X_W$.
	Set $\pi_W \colon X_W\to W$ to be the restriction of $\pi_{O_K}$ on $X_W$.
	Then $o_W$ is a section of $\pi_{O_K}$.
	For every $t\in W$, denote by $X_t$, $f_t$ and $o_t$ the specializations of $X_W$, $f_W$ and $o_W$ at $t$.
	After enlarging $B(f,\mathbb{Z})$,
	$X_t$ is irreducible for every $t\in W$.
	Then for every point $x\in X^{\reg}(\overline{\mathbb{Q}})$, if $f^m(x)=x$ for some $m\geq 1$ and $\beta_1,\beta_2$ are the eigenvalues of the tangent map $d(f^m)|_x$, then for every prime $p\not\in B(f,\mathbb{Z})$, and every embedding $\tau \colon \overline{\mathbb{Q}} \hookrightarrow \mathbb{C}_p$, we have $|\tau(\beta_1)|,|\tau(\beta_2)|\leq 1$.

	Since $f$ is amplified, $\deg(f|_{\overline{C}})\geq 2$ by \cite{xie2022existence}*{Lemma 5.2}.
	Then $\overline{C}$ is either $\mathbb{P}^1$ or an elliptic curve.
	Since on a complex elliptic curve, an endomorphism of degree at least $2$ is everywhere repelling, $\overline{C}$ could not be an elliptic curve.
	Then we have $\overline{C} \cong \mathbb{P}^1$.
	Since $0 < |\alpha(\mu)| < 1$, by \cite{milnor2006dynamics}*{Corollary 14.5}, $f|_{\overline C}$ is not post-critically finite.

	Denote by $J(f)$ the union of the critical locus of $f$ and the singular locus of $X$.
	Since $o \notin J(f)$ and $o\in C$, we have $C \not\subset J(f)$.
	Then $C \cap J(f)$ is finite.
	Let $P(f,C)$ be the union of the orbits of all periodic points in $C\cap J(f)$.
	Then $P(f,C)$ is finite.
	Observe that for every $n \geq 1$, $P(f^n,C) = P(f,C)$.
	By \cite{xie2022existence}*{Lemma 6.8},
	after replacing $f$ by a suitable positive iteration,
	there is a prime $p \notin B(f,\mathbb{Z})$,
	an embedding $\tau \colon \overline{\mathbb{Q}} \hookrightarrow \mathbb{C}_p$
	and some $\overline{x} \in \mathrm{Fix}(f|_{\overline{C}}) \setminus \pi_C^{-1}(P(f,C))$
	such that $C$ is smooth at $\pi_C(\overline{x})$
	and $|\tau(d(f|_{\overline{C}})|_{\overline x})| < 1$.
	Set $x \coloneqq \pi_C(\overline{x})$.
	Since $x \notin P(f,C)$, $X$ is smooth at $x$ and $df|_x$ is invertible.
	Since $d(f|_{\overline C})|_{\overline x}$ is an eigenvalue of $df|_x$,
	$x$ is a good fixed point of $f$.
\end{proof}

\begin{lemma}\label{lemextgoodfix}
	Assume $f$ is an amplified endomorphism which has R-property.
	Then either $(X,f)$ satisfies SAZD-property,
	or there is an $n\geq 1$ such that $f^n$ has a good fixed point.
\end{lemma}

\begin{proof}
	The proof is a small modification of \cite{xie2022existence}*{Lemma 6.5}.
	Indeed, replacing \cite{xie2022existence}*{Lemma 6.7} by \cref{prop_gd_pt},
	the proof of \cite{xie2022existence}*{Lemma 6.5} works.
\end{proof}

The following result is a singular version of \cite{xie2022existence}*{Proposition 6.15}.

\begin{proposition}\label{prozaridenseorbitsurfendoamp}
	Let $X$ be a projective surface over $\bk$, and $f \colon X \to X$ an amplified endomorphism.
	Assume that $f$ satisfies R-property.
	Then
	$(X,f)$ satisfies SAZD-property.
\end{proposition}

\begin{proof}
	By \cref{lemextgoodfix} and \cref{prop:AZD_bir}, we may assume that $f$ has a good fixed point.
	Then \cite{xie2022existence}*{\S 6.3},
	whose proof still works in the singular case,
	shows that the pair $(X,f)$ satisfies SAZD-property.
\end{proof}


Now we are ready to give the following:

\begin{proof}[Proof of \cref{prop:int_amp}]
	Pick any embedding $\sigma \colon \bk\hookrightarrow \mathbb{C}$.
	View $X_{\bk}(\mathbb{C})$ as a complex surface induced by $\sigma$.
	Let $\pi \colon X'\to X$ be a projective desingularisation of $X$.
	Set $f' \coloneqq \pi^{-1}\circ f\circ \pi \colon X'\dashrightarrow X'$.
	We have $\delta_{f'}=\delta_f$ and $\deg f'=\deg f$.
	Let $U$ be a Zariski open subset of $X'$ such that $\pi|_U$ is an isomorphism to its image.
	By \cite{guedj2005ergodic}*{Theorem~3.1, (iv)}, \cite{dinh2015equidistribution}*{Theorem~1.1}
	and since $\deg f>\delta_f$,
	there is an $m \geq 1$ and a repelling fixed point $o$ of $f'^m$ in $U$.
	Then $\pi(o)$ is a smooth repelling fixed point of $f^m$.
	Thus $f^m$ has R-property.
	By \cref{prop:AZD_bir,prozaridenseorbitsurfendoamp},
	$(X, f^m)$ and hence $(X,f)$ satisfy SAZD-property.
	The final assertion follows from \cref{rem:pol_to_amp}.
\end{proof}

The following is borrowed from \cite{xie2022existence}*{\S 7}.

\begin{proposition}\label{prop:AZD_ruled}
	Suppose that $\pi \colon X \to Y$ is a $\mathbb{P}^1$-bundle over a smooth projective curve $Y$ and a non-isomorphic surjective endomorphism $f \colon X \to X$ descends to an endomorphism $g \colon Y \to Y$.
	Then $(X, f)$ satisfies AZD-property.
\end{proposition}

\begin{proof}
	This is proved in \cite{xie2022existence}*{\S 7}.
	We sketch it here for the convenience of the reader.
	Write $\NE(X) = \langle [F], [E]\rangle$, where $F$ is a general fibre of $\pi$.
	Then $f^*F \equiv \delta_g F$.
	Write $f^*E \equiv \lambda E$.
	Then $\delta_f = \max\{\lambda, \delta_g\}$, and $1 < \deg f = \lambda \delta_g$.
	If both $\lambda > 1$ and $\delta_g > 1$, then $f$ is int-amplified so $(X, f)$ satisfies AZD-property by \cref{prop:int_amp}.

	Thus we may assume $\lambda = 1$ or $\delta_g = 1$.
	In particular, $\lambda \neq \delta_g$.
	Then $\lambda^2 E^2 = (f^*E)^2 = (\deg f) E^2 = (\lambda \delta_g) E^2$ implies that $E^2 = 0$.

	By \cref{prop:AZD_bir,prop:SAZD_inf_curve}
	and replacing $f$ by an iteration,
	we may assume that there are infinitely many (irreducible) curves $C_i$ with $f(C_i) = C_i$.
	For $C = C_i$, write $C = aF + bE$.
	Then $a \lambda F + b \delta_g E = f_* C$ (which is proportional to $C$)
	and $\lambda \ne \delta_g$ imply that $C = C_i$ is proportional to $F$ or $E$.

	Suppose that infinitely many $C_i$'s are proportional to $F$.
	Then $C_i \cdot F = 0$ and hence $C_i$ equals $X_{y_i}$, a fibre of $\pi$ over $y_i \in Y$.
	Now $f(C_i) = C_i$ implies that $g(y_i) = y_i$.
	Then $g$ has infinitely many fixed points $y_i$'s, so $g = \id_Y$.
	Hence $(X, f)$ satisfies AZD-property.

	Thus we may assume all $C_i$'s are proportional to $E$.
	Since $E^2 = 0$, all $C_i$'s are disjoint.
	Taking base changes by (the normalisation of) $C_i \to Y$ consecutively
	and by \cref{lem:AZD_gen},
	we may assume that $C_i$ ($i = 1, 2, 3$) are cross-sections of $\pi$.
	Then there is a natural isomorphism $Y \times \mathbb{P}^1 \to X$
	mapping $Y \times \{0, 1, \infty\}$ to $C_1 \cup C_2 \cup C_3$.
	Identifying $X = Y \times \mathbb{P}^1$,
	our $f = g \times h$ with $h \colon \mathbb{P}^1 \to \mathbb{P}^1$ a morphism.
	Now the still disjointness of all $C_i$ ($i \ge 2$) with $C_1$
	implies that $C_i$ equals $X_{p_i}$,
	a fibre of the projection $X \to \mathbb{P}^1$ over $p_i \in \mathbb{P}^1$.
	Also $f(C_i) = C_i$ implies $h(p_i) = p_i$ for infinitely many points $p_i$'s.
	Hence $h = \id_{\mathbb{P}^1}$.
	Thus $(X, f)$ satisfies AZD-property.
\end{proof}


\begin{proof}[Proof of \cref{thm:AZD}]
	By \cref{cor:AZD_aut} and \cref{prop:AZD_bir},
	we may assume $\deg f \ge 2$ and $X$ is a normal projective surface.
	By \cref{cor:main_compact}, \cref{prop:int_amp} and \cref{prop:AZD_bir}, \cref{lem:AZD_gen},
	we may assume either $X$ is an abelian surface,
	or $X$ is a $\mathbb{P}^1$-bundle and $f$ descends to the base.
	Then the theorem follows from \cite{xie2022existence}*{Theorem 1.14} and \cref{prop:AZD_ruled}.
\end{proof}

\begin{proof}[Proof of \cref{pro:AZDinforbits}]
	There is a finitely generated field extension $K$ over $\mathbb{Q}$
	such that $\overline{K}=\bk$, and $X, f$ are defined over $K$.
	There exists a subring $R$ of $K$ such that $R$ is finitely generated over $\mathbb{Z}$
	and $\Frac R=K$.
	Pick a model $\pi\colon X_R \to \Spec R$ which is projective over $\Spec R$
	and whose generic fibre is $X_K$.

	Our $f$ extends to a rational self-map $f_R\colon X_R\dashrightarrow X_R$.
	Denote by $B_R$ the indeterminacy locus of $f_R$.
	By \cite{xie2022existence}*{Lemma 3.23},
	there exists a nonempty, affine open subset $U$ of $\Spec R$ such that
	\begin{enumerate}
		\item $U$ is of finite type over $\Spec \mathbb{Z}$;
		\item for every point $y \in U$, the fibre $X_y$ is geometrically irreducible
		      and $\dim_{K(y)} X_y = \dim_K X_K$, where $K(y)$ is the residue field at $y$; and
		\item for every $y \in U$, the fibre $X_y$ is not contained in $B_{R}$
		      and the restriction $f_y$ of $f_R$ to $X_y$ is dominant.
	\end{enumerate}
	Moreover, after shrinking $U$, we may assume that for every $y\in U$, $f_y$ is separable.
	By \cite{bell2016dynamical}*{Proposition 2.5.3.1},
	there are infinitely many primes $p \geq 3$
	such that $R$ can be embedded into $\mathbb{Z}_p\subseteq \mathbb{C}_p^{\circ}$.
	This induces an embedding $\Spec \mathbb{Z}_p \to \Spec R$.
	Denote by $\tau\colon K\hookrightarrow \mathbb{C}_p$ the field embedding.
	Set $X_{\mathbb{C}_p^{\circ}} \coloneqq X_R \times_{\Spec R}\Spec \mathbb{C}_p^{\circ}$,
	and $f_{\mathbb{C}_p^{\circ}}\coloneqq f_R\times_{\Spec R} \id$.
	Let $(X_{\mathbb{C}_p}, f_{\mathbb{C}_p})$ and $(X_{\mathbb{F}_p}, f_{\mathbb{F}_p})$ be the generic fibre
	and special fibre of $(X_{\mathbb{C}_p^{\circ}}, f_{\mathbb{C}_p^{\circ}})$.
	Then $X_{\overline{\mathbb{F}_p}}$ is irreducible and $f_{\mathbb{F}_p}$ is dominant.
	Denote by $B_{\mathbb{C}_p^{\circ}}$ the base change of $B_R$.
	Then $X_{\overline{\mathbb{F}_p}}\not\subseteq B_{\mathbb{C}_p^{\circ}}$.

	Since $(X, f)$ has SAZD-property,
	there is a nonempty adelic open subset $A\subseteq X(\bk)$,
	such that for every $y\in A$, the orbit of $y$ is well-defined and Zariski dense.
	We need:

	\begin{lemma}
		\label{lemorbitfinitered}
		Let $K'/K$ be a finite extension and $\tau'\colon K'\hookrightarrow \mathbb{C}_p$ a field embedding extending $\tau$.
		Let $V$ be any nonempty Zariski open subset of $X_{\overline{\mathbb{F}_p}}\setminus B_{\mathbb{C}_p^{\circ}}$.
		Then there is $x\in A$ and a field embedding $\overline{\tau'}\colon \bk\hookrightarrow \mathbb{C}_p$ extending $\tau'$,
		such that the reduction of $\phi_{\overline{\tau'}}(O_f(x))$ to $X_{\overline{\mathbb{F}_p}}$ is finite and contained in $V$.
		Here $\phi_{\overline{\tau'}}\colon X(\bk)\hookrightarrow X(\mathbb{C}_p)$ is the embedding induced by $\overline{\tau'}$.
	\end{lemma}

	Assuming \cref{lemorbitfinitered},
	we first construct points $x_n\in A$ ($n\geq 1$),
	increasing finite extensions $K_n$ ($n\geq 1$) of $K$ over which $x_n$ is defined,
	and field embeddings $\overline{\tau_n}\colon \bk\hookrightarrow \mathbb{C}_p$
	with $\overline{\tau_n}|_{K_{n-1}}=\overline{\tau_{n-1}}|_{K_{n-1}}$,
	such that the reductions of $\phi_{\overline{\tau_n}}(O_f(x_n))$ to $X_{\overline{\mathbb{F}_p}}$
	are finite,
	contained in $X_{\overline{\mathbb{F}_p}}\setminus B_{\mathbb{C}_p^{\circ}}$ and disjoint.
	For these $x_n$,
	we have $\phi_{\overline{\tau_n}}(O_f(x_n))=\phi_{\overline{\tau_m}}(O_f(x_n))$ ($m\geq n$).
	Hence the orbits of $x_n$ ($n\geq 1$) are disjoint,
	thus, \textit{it proves \cref{pro:AZDinforbits}}.

	We construct $x_1, \overline{\tau_1}$ by applying \cref{lemorbitfinitered} to $V=X_{\overline{\mathbb{F}_p}}\setminus B_{\mathbb{C}_p^{\circ}}$ and $K'=K$.
	Let $K_1$ be any finite field extension of $K$ such that $x_1$ is defined over $K_1$.
	Assume that we have constructed $x_n, K_n, \overline{\tau_n}$ ($n=1,\dots, m$).
	Let $S_m$ be the union of the reductions of $\phi_{\overline{\tau_n}}(O_f(x_n))=\phi_{\overline{\tau_m}}(O_f(x_n))$
	($n=1,\dots, m$), which is a finite subset of $X_{\overline{\mathbb{F}_p}}\setminus B_{\mathbb{C}_p^{\circ}}$.
	We construct $x_{m+1}, \overline{\tau_{m+1}}$ by applying \cref{lemorbitfinitered} to $V=X_{\overline{\mathbb{F}_p}}\setminus (B_{\mathbb{C}_p^{\circ}}\cup S_m)$ and $K'=K_{m}$.
	Let $K_{m+1}$ be any finite field extension of $K_m$ such that $x_{m+1}$ is defined over $K_{m+1}$.
\end{proof}

\begin{proof}[Proof of \cref{lemorbitfinitered}]
	Applying \cite{amerik2011existence}*{Corollary 2} to the rational self-map $f_{\mathbb{F}_p}|_{V}\colon V\dashrightarrow V$, there exists a periodic point $\bar{x}\in V$ whose orbit under $f_{\mathbb{F}_p}$ is contained in $V$.
	Let $U$ be the $p$-adic open subset of $X(\mathbb{C}_p)$ of points whose reduction is $\overline{x}$.
	Let $X_{K'}(\tau', U)$ be the basic adelic subset over $K'$ associated to $\tau'$ and $U$ as defined in \cite{xie2022existence}*{Section 3.11}.
	It is a nonempty adelic open subset of $X(\bk)$.

	Since $X$ is irreducible, $A\cap X_{K'}(\tau', U)\neq\emptyset$.
	Pick $x\in A\cap X_{K'}(\tau', U)$.
	Then the orbit of $x$ is well-defined and Zariski dense.
	By definition of $X_{K'}(\tau', U)$, some field embedding $\overline{\tau'}\colon \bk\hookrightarrow \mathbb{C}_p$ extends $\tau'$ with $\phi_{\overline{\tau'}}(x)\in U$.
	Since the reduction of $\phi_{\overline{\tau'}}(x)$ to $X_{\overline{\mathbb{F}_p}}$ is $\overline{x}$ and the orbit of $\overline{x}$ is finite and contained in $V$, this proves \cref{lemorbitfinitered} and also \cref{pro:AZDinforbits}.
\end{proof}

\end{document}